\documentclass{article}

\usepackage{here}
\usepackage[utf8]{inputenc}
\usepackage[T1]{fontenc}
\usepackage[british]{babel}
\usepackage{enumitem}
\usepackage{authblk}

\usepackage{amsmath}
\usepackage{graphicx}
\usepackage{subcaption}
\usepackage{dsfont}
\usepackage[marginal]{footmisc}
\usepackage[all,cmtip]{xy}
\usepackage{amssymb}
\usepackage{caption}
\usepackage{xcolor}
\usepackage{amsthm}
\theoremstyle{definition}
\newtheorem{defn}{Definition}[section]

\theoremstyle{plain}
\newtheorem{thm}[defn]{Theorem}

\newtheorem{lem}[defn]{Lemma}

\newtheorem{prop}[defn]{Proposition}

\theoremstyle{remark}

\newtheorem{rem}[defn]{Remark}

\theoremstyle{definition}

\DeclareMathOperator{\dist}{dist}
\DeclareMathOperator{\vol}{vol}

\DeclareMathOperator{\var}{var}

\DeclareMathOperator{\inn}{int}

\DeclareMathOperator{\Res}{Res}
\DeclareMathOperator{\spann}{span}
\usepackage{mathrsfs}

\usepackage{tikz-cd}
\usepackage{tikz}
\usetikzlibrary{arrows,topaths,calc,patterns,angles,quotes}

\usetikzlibrary{decorations.fractals}
\usetikzlibrary{arrows,calc,patterns}

\newcommand{\e}{\epsilon}

\renewcommand{\epsilon}{\varepsilon}

\newcommand{\gen}{G}
\newcommand{\coeffweyl}{Q}
\newcommand{\coeffmink}{R}

\usepackage{hyperref}

\numberwithin{equation}{section}

\allowdisplaybreaks
\usepackage[left=3cm,right=3cm,bottom=2cm,top=2cm]{geometry}
\begin{document}
\title{Eigenvalue counting functions and parallel volumes for examples of fractal sprays generated by the Koch snowflake}
\author[1]{S. Kombrink\footnote{\url{s.kombrink@bham.ac.uk}}}
\author[1]{L. Schmidt\footnote{\url{l.schmidt@bham.ac.uk}. 
The second author was supported by EPSRC DTP and the University of Birmingham.\par~\, The authors would like to thank the anonymous referee for helpful comments.}}

\affil[1]{School of Mathematics, University of Birmingham, Edgbaston, Birmingham, B15\,2TT, UK}
\makeatletter
\newcommand{\subjclass}[2][2020]{%
  \let\@oldtitle\@title%
  \gdef\@title{\@oldtitle\footnotetext{#1 \emph{Mathematics subject classification.} #2}}%
}
\newcommand{\keywords}[1]{%
  \let\@@oldtitle\@title%
  \gdef\@title{\@@oldtitle\footnotetext{\emph{Key words and phrases.} #1.}}%
}
\makeatother
\date{}
\maketitle
\vspace{-0.5cm}
\begin{abstract}
    We apply recent results by the authors to obtain bounds on remainder terms of the Dirichlet Laplace eigenvalue counting function for domains that can be realised as countable disjoint unions of scaled Koch snowflakes. Moreover we compare the resulting exponents to the exponents in the asymptotic expansion of the domain's inner parallel volume.
\end{abstract}

\section{Overview}
Let $\Omega$ be an open, bounded subset of $\mathbb{R}^n$ and let $\partial \Omega$ denote its boundary. Defining the Laplace operator $\Delta:=\sum_i \partial_i ^2$ on $\Omega$, one may study the classical Laplace eigenvalue problem $-\Delta u = \lambda u$ in $\Omega$ under different boundary conditions, such as Dirichlet ($u = 0$ on $\partial\Omega$) or Neumann ( $\frac{\partial u}{\partial \mathbf{n}} = 0$ on $\partial\Omega$, where $\mathbf{n}$ denotes the exterior normal to $\partial \Omega$). Starting with Weyl's famous work on the asymptotic distribution of eigenvalues \cite{weyl1911}, much progress has been made concerning the relation between the geometry of $\Omega$  
and the spectrum of the Laplace operator
on $\Omega$ (i.\,e.\ the solutions $\lambda$ to the above Laplace eigenvalue problem)
under the condition that $\partial \Omega$ is sufficiently regular; see for example \cite{ivrii1,ivrii2,sogge1993}. For sufficiently regular $\partial \Omega$ the number (with multiplicity) of eigenvalues $\leq t$ subject to Dirichlet boundary conditions, $N_D(\Omega,t)$, is asymptotically given by
\begin{align}\label{eq:Weyl}
    N_D(\Omega,t) = (2\pi)^{-n} V_n \vol(\Omega)t^{n/2} - \frac{1}{4}(2\pi)^{-(n-1)} V_{n-1} \vol_{n-1}(\partial \Omega) t^{(n-1)/2} + o\left( t^{(n-1)/2} \right),
\end{align}
as $t\to\infty$, where $V_n$ is the volume of the $n$-dimensional unit ball. Originating with Berry's conjecture on a generalisation of \eqref{eq:Weyl} 
to domains with rough boundary in \cite{berry1979,berry1980}, significant progress has been achieved concerning domains with fractal boundary (see for example \cite{la1991,netrusov2007,NetrusovSafarov2005} and more recently \cite{teplyaev2022,hinzpierratteplyaev2023}). This quickly lead to the \emph{Weyl-Berry conjecture}, see \cite{lapidusdundee}. One variant of this conjecture states that whenever $\partial \Omega$ has Minkowski dimension $\dim_M(\partial\Omega)=d$ and is Minkowski measurable (i.\,e.\ its $d$-Minkowski content $M$ exists and is positive and finite), then the counting function satisfies
\begin{align}\label{eq:Weyl-Berry}
    N_D(\Omega,t) = (2\pi)^{-n} V_n \vol_n(\Omega)t^{n/2} + c_{n,d}M t^{d/2} + o(t^{d/2}),
\end{align}
as $t \to \infty$, where $c_{n,d}$ only depends on $n$ and $d$. A similar version of the conjecture has been formulated for the case of Neumann boundary conditions. 
Lapidus and Pomerance verified the conjecture \eqref{eq:Weyl-Berry} for $n=1$ (see \cite{lapo1993} or \cite{lapidusfrankenhuijsen2006} and the references therein as well as \cite{kenneth}). However, this conjecture is known to be incorrect in higher dimensions (see \cite{lapo1996} and \cite{huasleemann}).\par
In many cases a domain $\Omega$ with self-similar boundary can be understood as a \emph{fractal spray}, i.\,e.\ a disjoint union of rescaled copies of some fundamental domain (the \emph{generator}). In regards to the Laplace eigenvalue problem, generators are typically assumed to be of very regular nature. Aspects of renewal theory appear naturally in this context, see for example \cite{huasleemann98,kombrinkkesse2017,lapidusfrankenhuijsen2006}.
Notable examples of generators used in the literature are intervals (for fractal strings), squares (for carpets such as the Sierpi\'{n}ski carpet) or triangles (for gaskets such as the Sierpi\'{n}ski gasket). \par
In the present article we extend the investigations to generators with irregular boundary. 
More precisely, here, we consider a family of lattice fractal sprays $\Omega(k_1,k_2)$ in $\mathbb R^2$ generated by the Koch snowflake $K$ with $\dim \partial K=\delta = \log_3 4$ (see Sec.~\ref{sec:kochspray}, Fig.~\subref{fig:kochflockesec5} and Fig.~\subref{fig:kochflockesec5k1k2}) through a set of similarity maps $\{\phi_i\}_{i \in \Sigma}$ for which there exists $a>0$ and $\nu_i \in \mathbb{Z}$  s.\,t.\ $r_i:=|\phi_i '| = {\textup{e}}^{-a \nu_i}$. If $a>0$ is maximal with the property that the scaling ratios $r_i$ are all powers of ${\textup{e}}^{-a}$ then $a$ is called the \emph{lattice constant} of $\{\phi_i\}_{i \in \Sigma}$.  For such lattice fractal sprays $\Omega(k_1,k_2)$ we give an asymptotic expansion with error term of $N_D(\Omega(k_1,k_2),t)$ in Sec.~\ref{sec:weylasymp} and an asymptotic expansion of the inner parallel volume in Sec.~\ref{sec:minkcontent}. Both results are based on similar ideas from renewal theory. The exponents in the expansions correspond to poles of transfer operators, namely to certain elements of the sets
\begin{align*}
 \mathcal{Z}_{\text{C}} &:=\bigg{\{} z\in\mathbb C \ :\  \sum_{i \in \Sigma} r_i^{-2z} =1,\ \Im(z)\in\big{[}0,\frac{\pi}{a}\big{)} \bigg{\}} \qquad \text{ and }\\
 \mathcal{Z}_{\text{P}} &:=\bigg{\{} z\in\mathbb C \ :\  \sum_{i \in \Sigma} r_i^{2-z} =1,\ \Im(z)\in\big{[}0,\frac{2\pi}{a}\big{)} \bigg{\}}
\end{align*}
respectively.
The resulting expansion for the counting function is given by
\begin{align*}
    &N_D(\Omega,t) = \frac{1}{4\pi} \vol_2(\Omega(k_1,k_2)) t
    - \sum_{ \substack{z \in \mathcal{Z}_{\text{C}}:\\\Re (z) < -\delta/2}} A_{z,\widetilde{\beta}}\cdot t^{-z} + {o}\left( t^{\delta/2+\gamma} \right),
\end{align*}
as $t\to\infty$ for any $\gamma>0$,
where $\lvert A_{z,\widetilde{\beta}}\rvert$ are bounded and oscillatory in $\widetilde{\beta} = 2a\{ \log t/2a \}$, with $\{x\}$ denoting the fractional part of $x\in\mathbb R$. We denote the inner $\epsilon$-parallel set by $\Omega_{-\epsilon}:=\{x \in \Omega:\dist(x,\partial \Omega)<\epsilon\}$. For the inner $\epsilon$-parallel volume the expansion is
\begin{align*}
    &\vol_2( \Omega_{-\epsilon}) 
     = B_{\widetilde{\beta},\left(2- \delta\right)}\epsilon^{2-\delta}
     +B_{\widetilde{\beta},2}\epsilon^{2}  
     +\sum_{z\in \mathcal{Z}_{\text{P}}} B_{\widetilde{\beta},z}\epsilon^{z}+o(\epsilon^{\gamma})
\end{align*}
as $\e\to 0$ for any $\gamma>0$, where $B_{\widetilde{\beta},z}$ are bounded and oscillatory in $\widetilde{\beta}=a\{ -\log\epsilon/a \}$. 
Remarkably, the one-to-one correspondence between elements in $\mathcal Z_{\text{C}}$ and in $\mathcal Z_{\text{P}}$ shows a strong connection between the eigenvalue counting function on $\Omega(k_1,k_2)$ and the inner $\epsilon$-parallel volume of $\Omega(k_1,k_2)$.

\section{Fractal sprays generated by the Koch snowflake}\label{sec:kochspray}
\begin{figure}[htb]
    \begin{minipage}[t]{.45\textwidth}
        \centering
            \begin{tikzpicture}[scale=3.5,decoration=Koch snowflake]
    \draw[scale={sqrt(3)},shift={({-0.5/sqrt(3)},{-1/6})},rotate=30] decorate{decorate{decorate{decorate{decorate{ (0,0) -- (1,0) -- (0.5,{-sqrt(3)/2}) -- (0,0)}}}}};
    \draw decorate{decorate{decorate{decorate{decorate{ (0,0) -- (1,0) -- (0.5,{-sqrt(3)/2}) -- (0,0)}}}}};
    \draw[fill=gray,scale={sqrt(3)/3},shift={({sqrt(3)/6},-0.5)},rotate=30] decorate{decorate{decorate{decorate{decorate{ (0,0) -- (1,0) -- (0.5,{-sqrt(3)/2}) -- (0,0)}}}}};
    \node at (0.5,-0.27) {$K$};
    \node at (0.17,0.3) {$\phi_{3}$};
    \node at (0.83,0.3) {$\phi_{2}$};
    \node at (1.17,-0.28) {$\phi_{1}$};
    \node at (0.83,-0.85) {$\phi_{6}$};
    \node at (0.17,-0.85) {$\phi_{5}$};
    \node at (-0.17,-0.28) {$\phi_{4}$};
    \node at (0.5,0.1) {\small{$\phi_{8}$}};
    \node at (0.83,-0.075) {\small{$\phi_{7}$}};
    \node at (0.83,-0.47) {\small{$\phi_{12}$}};
    \node at (0.5,-0.65) {\small{$\phi_{11}$}};
    \node at (0.17,-0.47) {\small{$\phi_{10}$}};
    \node at (0.17,-0.075) {\small{$\phi_{9}$}};
    \draw[<->] (1.55,{0}) -- (1.55,{-7/12});
    \node[right] at (1.55,{0.5*(3.5/6-7/6)}) {$b$};
    \end{tikzpicture}
        \subcaption{Depiction of the IFS of the fractal spray studied in Sec.~\ref{sec:weylasymp} and Sec.~\ref{sec:minkcontent}. The \emph{base length} $b$ of the snowflake $K$ is also shown.}
      \label{fig:kochflockesec5}
    \end{minipage}
    \hfill
    \begin{minipage}[t]{.45\textwidth}
        \centering
            \begin{tikzpicture}[scale=3.5,decoration=Koch snowflake]
    \draw[scale={sqrt(3)},shift={({-0.5/sqrt(3)},{-1/6})},rotate=30] decorate{decorate{decorate{decorate{decorate{ (0,0) -- (1,0) -- (0.5,{-sqrt(3)/2}) -- (0,0)}}}}};
    \draw decorate{decorate{decorate{decorate{decorate{ (0,0) -- (1,0) -- (0.5,{-sqrt(3)/2}) -- (0,0)}}}}};
    \draw[fill=gray,scale={sqrt(3)/3},shift={({sqrt(3)/6},-0.5)},rotate=30] decorate{decorate{decorate{decorate{decorate{ (0,0) -- (1,0) -- (0.5,{-sqrt(3)/2}) -- (0,0)}}}}};
    \node at (0.5,-0.29) {$K_0$};
    \node at (0.17,0.3) {$\phi_{3}$};
    \node at (0.83,0.3) {$\phi_{2}$};
    \node at (0.83,-0.85) {$\phi_{6}$};
    \node at (0.17,-0.85) {$\phi_{5}$};
    \node at (-0.17,-0.28) {$\phi_{4}$};
    \node at (0.5,0.1) {\small{$\phi_{8}$}};
    \node at (0.83,-0.075) {\small{$\phi_{7}$}};
    \node at (0.83,-0.47) {\small{$\phi_{12}$}};
    \node at (0.5,-0.65) {\small{$\phi_{11}$}};
    \node at (0.17,-0.47) {\small{$\phi_{10}$}};
    \node at (0.17,-0.075) {\small{$\phi_{9}$}};
    \draw[fill=gray!50!white,scale={1/3},shift={({3},{-sqrt(3)/3})},rotate=0] decorate{decorate{decorate{decorate{decorate{ (0,0) -- (1,0) -- (0.5,{-sqrt(3)/2}) -- (0,0)}}}}};
    \node at ({0.98+sqrt(3)/9},-0.29) {$K_1$};
    \node at (1.4,-0.28) {\tiny$\phi_{1,1}$};
    \node at (1.27,-0.1) {\tiny$\phi_{1,2}$};
    \node at (1.05,-0.1) {\tiny$\phi_{1,3}$};
    \node at (0.95,-0.28) {\tiny$\phi_{1,4}$};
    \node at (1.05,-0.48) {\tiny$\phi_{1,5}$};
    \node at (1.27,-0.48) {\tiny$\phi_{1,6}$}; 
    \end{tikzpicture}
        \subcaption{Depiction of the variant $\Omega(1,0)$ of the IFS of the fractal spray studied in Sec.~\ref{sec:weylasymp} and Sec.~\ref{sec:minkcontent}. In this case, the map $\phi_1$ has been replaced by the six maps $\phi_{1,1},\ldots,\phi_{1,6}$ giving rise to an additional connected component of the generator $\gen=K_0 \cup K_1$. Analogously one can replace $\phi_7,\ldots,\phi_{12}$.}
      \label{fig:kochflockesec5k1k2}
    \end{minipage}  
    \label{fig:1-2}
\end{figure}
\begin{figure}
\centering
 \begin{subfigure}[t]{0.32\textwidth}
    \centering
    \includegraphics[width=\textwidth]{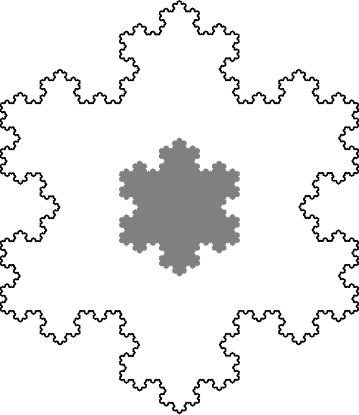}
 \end{subfigure}
 \hfill
 \begin{subfigure}[t]{0.32\textwidth}
 \centering
     \includegraphics[width=\textwidth]{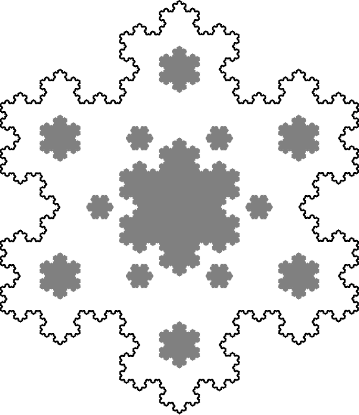}
 \end{subfigure}
 \hfill
 \begin{subfigure}[t]{0.32\textwidth}
 \centering
     \includegraphics[width=\textwidth]{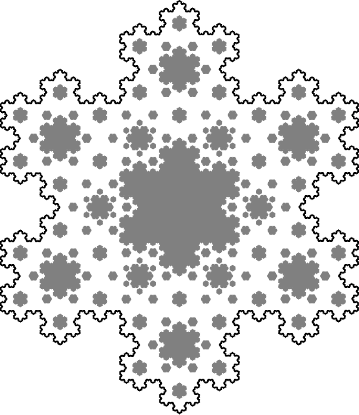}
 \end{subfigure}
    \caption{Example of the iterative construction of $\Omega$ as defined in \eqref{eq:Omega}. From left to right the $0^{\text{th}}$, first and second iterations of the generator $K$ under $\Phi$ are shown, rotated by $30^\circ$. The $0^{\text{th}}$ iteration (left) shows $K$. The first iteration (middle) shows $K \cup \Phi(K)$  . The second iteration (right) then shows 
    $K \cup \Phi(K) \cup \Phi^2(K)$.} 
    \label{fig:fractalexample}
\end{figure}

  
  \noindent{}Our class of examples is based on the construction shown in Fig.~\subref{fig:kochflockesec5} (see also Fig.~\ref{fig:fractalexample}). More precisely, we consider the iterated function system (IFS) $\Phi:=\Phi(0,0):=\{\phi_1,\ldots,\phi_{12}\}$ defined on $\mathbb R^2$ given by the maps
\begin{align*}
    \phi_i(x) := 
    \begin{cases}
    \frac{1}{3}x + \frac{2\sqrt{3}}{3}\begin{pmatrix}
        \cos [(i-1)\pi/3]\\
        \sin [(i-1)\pi/3]
    \end{pmatrix} &\text{ if }i \in \{1,\ldots,6\}\\
    \frac{1}{3\sqrt{3}}R_{\pi/6}(x) + \frac{2}{3}\begin{pmatrix}
        \cos[\pi/6 + (i-7)\pi/3]\\
        \sin[\pi/6 + (i-7)\pi/3]
    \end{pmatrix} &\text{ if }i \in \{7,\ldots,12\}
     \end{cases}
\end{align*}
where $R_{\pi/6}$ is a rotation by $\pi/6$ about the origin. Define the action of $\Phi$ on subsets of $\mathbb R^2$ by $\Phi A:=\bigcup_{i\in\Sigma} \phi_i A$ with $i\in \{1,\ldots,12\}=:\Sigma$. Further, let $F$ denote the unique non-empty compact invariant set associated to the IFS, i.\,e.\ the set satisfying $F=\Phi F := \bigcup_{i\in\Sigma}\phi_i F$. 
Note that $F$ is contained in the disk around the origin of radius $\sqrt{3}$, with $\binom{-\sqrt{3}}{0}$ and $\binom{\sqrt{3}}{0}$ belonging to $F$.
As all $\phi_i$ are similarities, $F$ is a self-similar set. From the definition of the maps it is evident that $\Phi$ consists of six contractions with contraction ratios $r_{1},\ldots,r_6 = 1/3=(\exp(-a))^2$ and six contractions with contraction ratios $r_{7},\ldots,r_{12} = \sqrt{3}/9=(\exp(-a))^3$, where $a:=\log{3}/2$ is known as the {\em lattice constant}.
Note that by construction, $\mathbb R^2\setminus F$ has got a unique unbounded connected component, which we denote by $U$. Further, we let $O:=\mathbb R^2 \setminus \overline{U}$. Then $O$ is open and satisfies $\phi_i O\subseteq O$ for $i\in\Sigma$ and $\phi_i O\cap \phi_j O=\emptyset$ for $i\neq j\in\Sigma$. This implies that the open set condition (OSC) is satisfied and that $O$ is a feasible open set for the OSC. Fig.~\subref{fig:kochflockesec5} shows the images of $\overline{O}$ under the maps $\phi_1,\ldots,\phi_{12}$. 

For $m\in\mathbb N$ and $\omega=(\omega_1,\ldots,\omega_m)\in \Sigma^m$, let $\phi_{\omega}:=\phi_{\omega_1}\circ\cdots\circ \phi_{\omega_m}$ and 
\begin{equation}\label{eq:nu}
    r_{\omega}
    :=\exp(-a\nu_{\omega}):= \prod_{i=1}^m r_i
    \quad\text{with}\quad \nu_{\omega}:=\sum_{i=1}^m \nu_{\omega_i}\in \mathbb N.
\end{equation}
We define $K:=O\setminus \overline{\Phi O}$, and note that $K$ is the interior of a Koch snowflake with base length $b=1$, see Fig.~\subref{fig:kochflockesec5}. A central object of our studies is 
\begin{equation}\label{eq:Omega}
    \Omega := \bigcup_{m=0}^{\infty} \Phi^m(K) 
    := K \cup \bigcup_{m=1}^{\infty} \bigcup_{\omega\in\Sigma^m} \phi_{\omega}(K)
\end{equation}
which is a countable union of disjoint open sets $\phi_w (K)$ and can be viewed as a fractal spray with \emph{generator} $K$. Here, $\Phi^0(K)$ is understood to be $K$. The first three iterations of the construction of $\Omega$ in \eqref{eq:Omega} are shown in Fig.~\ref{fig:fractalexample}.\par
We will moreover study the sets $\Omega(k_1,k_2)$ which result from modifications of the above construction as explained below. For $(k_1,k_2) \in \{0,\ldots,6\}^2$, we replace each of $\phi_1,\ldots,\phi_{k_1}$ with six maps of contraction ratio $1/9$ and each $\phi_7,\ldots,\phi_{k_2+6}$ with six maps with contraction ratio $1/(9\sqrt{3})$. The replacement of $\phi_i$ with six maps $\phi_{i,1},\ldots,\phi_{i,6}$ is done in such a way that $\bigcup_{k=1}^6 \phi_{i,k}O\subset \phi_i O$, that $\phi_{i,k} O \cap \phi_{i,j} O =\emptyset$ for all $k\neq j$, and that $\partial \phi_i O\subset \partial\bigcup_{k=1}^6 \phi_{i,k}O$. See Fig.~\subref{fig:kochflockesec5k1k2} for an example of the replacement procedure.
The corresponding IFS consisting of $12+5(k_1+k_2)$ maps will be denoted by $\Phi(k_1,k_2)$ and the associated alphabet by $\Sigma(k_1,k_2)$. The generator $O \setminus \overline{\Phi(k_1,k_2)(O)}$, that we denote by $\gen$ in this setting, has $1+k_1+k_2$ connected components. 
The fractal spray $\bigcup_{m=0}^{\infty}\Phi^m(k_1,k_2)(\gen)$ generated by $\gen$ 
will be denoted by $\Omega(k_1,k_2)$. 
We write $k_1=0$ when no replacement is intended for $\phi_1,\ldots,\phi_6$ (and  correspondingly $k_2=0$) so that $\Omega(0,0):=\Omega$.

 \section{Background on counting functions of the Koch snowflake}\label{sec:preliminaries} 
 Let $\Omega\subset\mathbb R^n$ be a bounded domain, i.\,e.\ a bounded open subset of $\mathbb{R}^n$ and denote its boundary by $\partial \Omega$. In case of Neumann boundary conditions we will need this domains to have finitely many connected components.
 By $H^1(\Omega)$ we denote the usual Sobolev space, i.\,e.\ the set of all $u \in L^2(\Omega)$ with a weak derivative $\nabla u \in L^2(\Omega)$. We equip $H^1(\Omega)$ with the usual inner product $(u,v)_{H^1(\Omega)} := (u,v)_{L^2(\Omega)} + (\nabla u,\nabla v)_{L^2(\Omega)}$ so that $H^1(\Omega)$ becomes a Hilbert space. Further, $H^1 _0 (\Omega)\subset H^1(\Omega)$ will denote the closure of the set of compactly supported $C^\infty(\Omega)$-functions in $H^1(\Omega)$. 
 On $H^1(\Omega)$, resp. $H^1 _0(\Omega)$ we consider the Laplacian $\Delta := \sum_{i=1} ^n \partial_i ^2$ and focus on the eigenvalue equation of $-\Delta$ subject to Neumann \eqref{eq:ewgl} or Dirichlet \eqref{eq:ewgld} boundary conditions.
 
 \begin{minipage}{0.4\textwidth}
\begin{align}
\begin{cases}
    -\Delta u = \lambda u &\text{ in } \Omega\\
    \frac{\partial u}{\partial \mathbf{n}} = 0 & \text{ on } \partial\Omega
\end{cases} \label{eq:ewgl}
\end{align}
 \end{minipage}
 \hspace{0.1\textwidth}
 \begin{minipage}{0.4\textwidth}
\begin{align}
\begin{cases}
    -\Delta u = \lambda u & \text{ in } \Omega\\
    u = 0 & \text{ on } \partial\Omega
\end{cases} \label{eq:ewgld}
\end{align}
 \end{minipage}~\\
\noindent{where $\mathbf{n}$ denotes the exterior normal to $\partial \Omega$.}
The variational formulation of the problem \eqref{eq:ewgl} is stated as follows: Find $u \in H^1(\Omega)$ s.\,t.\ $(\nabla u,\nabla v)_{L^2{(\Omega)}} = \lambda (u,v)_{L^2{(\Omega)}}$ for all $v \in H^1(\Omega)$. Note that the space in which this problem is studied dictates the boundary condition and that the variational formulation of \eqref{eq:ewgld} is: Find $u \in H^1_0(\Omega)$ s.\,t.\ $(\nabla u,\nabla v)_{L^2{(\Omega)}} = \lambda (u,v)_{L^2{(\Omega)}}$ for all $v \in H^1_0(\Omega)$. 
Replacing $H^1(\Omega)$ resp.\ $H^1_0(\Omega)$ with any other closed space $V$ satisfying $H^1 _0(\Omega) \subset V \subset H^1(\Omega)$ gives rise to variational problems with more general boundary conditions.
The corresponding (non-negative) spectrum will be denoted by $\sigma(-\Delta)$ and the essential spectrum by $\sigma_{\text{ess}}(-\Delta)$.\par
In order to define a counting function $N(\Omega,t) := \#\{\lambda \in \sigma(-\Delta):\lambda \leq t\}$, it is necessary that $\sigma(-\Delta)$ is discrete with the only accumulation point being at $\infty$. While this is always satisfied in case of Dirichlet boundary conditions, there are several occasions where this may fail to be true under Neumann boundary conditions as the essential spectrum can be non-empty in this case. Indeed, it was shown in \cite{hesesi1990} that any closed subset of $\mathbb{R}_{\geq 0}$ can be realised as the essential spectrum of a Laplacian on a domain $\subset \mathbb{R}^2$ subject to Neumann boundary conditions. On the other side, several criteria have been found which ensure that the essential spectrum is empty, see for example \cite{la1991,netrusov2007}. In particular, the Neumann Laplacian on domains bounded by quasicircles has vanishing essential spectrum. 
Consequently, as the Koch snowflake is a quasidisk, we  know that the essential spectrum of its Neumann Laplacian vanishes. In this context it is worth to mention the work of Rohde \cite{rohde2001} who showed that quasicircles are Rohde-snowflakes up to bi-Lipschitz transformations. Note that in this context, the unit square is understood to be a Rohde-snowflake.\par
In particular, for $V \in \{H^1 _0 (\Omega), H^1(\Omega)\}$ such eigenvalues are non-negative and by a variational argument the $k$-th eigenvalue of the problem \eqref{eq:ewgl}, resp. \eqref{eq:ewgld}, is given by
\begin{align}
    \lambda_k = \inf_{u \in V \cap \spann(u_1,u_2,\dots,u_{k'-1})^\perp} \frac{\|\nabla u\|^2}{\|u\|^2}, \label{eq:variation}
\end{align}
where $\{u_1,\dots,u_{k'-1}\}$ is an orthogonal basis of all eigenfunctions to eigenvalues $\lambda_1,\dots,\lambda_{k-1}$.
Additionally, there is a simple correspondence between counting functions and eigenvalues: Let $\Omega_1,\Omega_2$ be domains and $\lambda_k ^i$ be the $k$-th eigenvalue on $\Omega_i$. Then $(\lambda_k ^1 \leq \lambda_k ^2\,\, \forall k) \Leftrightarrow N(\Omega_1,t) \geq N(\Omega_2,t) \,\,\forall t$ and this is true for any considered boundary condition. We write $N_N(\Omega,t)$ (resp. $N_D(\Omega,t)$) for the number of eigenvalues (with multiplicity) of $-\Delta$ on $\Omega$ subject to Neumann (resp. Dirichlet) conditions on $\partial \Omega$.\par
One can deduce the following statements from \eqref{eq:variation}.
\begin{enumerate}
    \item The $k$-th Dirichlet eigenvalue $\lambda_k ^D$ and the $k$-th Neumann eigenvalue $\lambda_k ^N$ satisfy $\lambda_k ^D \geq \lambda_k ^N$. This is because $H^1 _0(\Omega) \subset H^1(\Omega)$ so that the infimum is taken over a larger set. In other words, $N_D(\Omega,t) \leq N_N(\Omega,t)$.\label{item:1}
    \item Let $\overline{\Omega} = \overline{\Omega_1} \cup \overline{\Omega_2}$ with $\Omega_1 \cap \Omega_2 = \emptyset$ and $\partial \Omega = \partial \Omega_1 \cup \partial \Omega_2$. Then
    \begin{align*}
        N_D(\Omega,t) = N_D(\Omega_1,t) + N_D(\Omega_2,t).
    \end{align*}
    With \eqref{eq:variation}, the reason lies in the existence of an isometric isomorphism $\iota:H^1 _0(\Omega_1 \sqcup \Omega_2) \simeq H^1 _0(\Omega_1) \oplus H^1 _0(\Omega_2)$ via $\iota:u \mapsto (u|_{\Omega_1} , u|_{\Omega_2})$ with inverse being $\iota^{-1}(u_1,u_2)(x) := u_i(x)$ whenever $x \in \Omega_i$.\label{item:2} 
    \item Writing $\alpha \Omega:= \{x \in \mathbb{R}^n: \alpha^{-1}x \in \Omega\}$ for $\alpha >0$, one has $N(\alpha \Omega,t) = N(\Omega,\alpha^{2}t)$ for both Dirichlet and Neumann boundary conditions. This is because of an isomorphism $\alpha^\ast\colon H^1(\Omega) \to H^1(\alpha\Omega),\ u\mapsto u\circ \alpha^{-1}$ and equally $\alpha^\ast\colon H^1 _0\Omega) \to H^1 _0(\alpha\Omega),\ u\mapsto u\circ \alpha^{-1}$.\label{item:3}
\end{enumerate}
Related to this is the Dirichlet-Neumann bracketing technique.

\paragraph{Dirichlet-Neumann bracketing.} 
  For a domain $\Omega \subset \mathbb{R}^n$ a \emph{volume cover} $\{\Omega_i\}_{i \in I}$ of $\Omega$ consists of at most countably many open sets $\Omega_i \subset \Omega$ with $\vol_n(\Omega) = \vol_n\left( \bigcup_{i \in I}\Omega_i \right)$.  Apart from the well-known classical Dirichlet-Neumann bracketing linking the counting functions of Dirichlet and Neumann eigenvalues, we mention a version that allows for non-disjoint covers as long as the elements of the cover do not intersect too often. More precisely one has the following result which also follows from the Min-Max-Principle. For any volume cover $\{\Omega_i\}_{i \in I}$ of $\Omega$, let $\mu:=\sup_{x \in \Omega}\#\{i \in I:x \in \Omega_i\}$ denote its \emph{multiplicity}.
  \begin{prop}[Dirichlet-Neumann bracketing with multiplicity, \cite{NetrusovSafarov2005}]\label{prop:bracketing}
   Let $\{\Omega_i\}_{i \in I}$ be a volume cover of $\Omega$. If the $\Omega_i$ are pairwise disjoint, then
   \begin{align*}
    \sum_{i \in I} N_D(\Omega_i,\lambda) \leq N_D(\Omega,\lambda) \leq N_N(\Omega,\lambda) \leq \sum_{i \in I} N_N(\Omega_i,\lambda)
   \end{align*}
   If the volume cover has finite multiplicity $\mu$, then
   \begin{align*}
    N_{N}(\Omega,\lambda) \leq \sum_{i \in I} N_{N}(\Omega_i,\mu \lambda).
   \end{align*}
  \end{prop}
 \subsection{Results for counting functions on the Koch snowflake}
  \begin{thm}[cf.~\cite{DOKUMENT}]\label{thm:satz1}
  Let $K$ be a Koch snowflake of base length $1$ as defined in Fig.~\subref{fig:kochflockesec5} and $\dim_H \partial K = \delta = \log_3 4$. Then
  \begin{align*}
     C_- \lambda^{\delta/2} \leq N_N(K,\lambda) -  \frac{1}{4\pi}\vol_2(K)\lambda \leq C_+ \lambda^{\delta/2}.
  \end{align*}
  for all $\lambda \geq 0.1$ with $C_-:= -1481$ and $C_+ :=281.5 \cdot 10^3$.
  \end{thm}
  \begin{proof}[Sketch of proof]
  The lower bound in the claim follows from a result in \cite{vandenBerg2001} since $N_D(\Omega,t) \leq N_N(\Omega,t)$. It actually holds true for all $\lambda \geq 0$.\par
  For the upper bound we follow an argument similar to \cite{NetrusovSafarov2005} and use the estimates on the upper inner Minkowski content found in \cite{LapidusPearse}. Let $\epsilon>0$ be sufficiently small and let $k$ be so that $\epsilon \in (3^{-(k+1)}/\sqrt{3},3^{-k}/\sqrt{3})$. We introduce a volume cover of $K_{-\epsilon}$ with domains $D_i ^\epsilon$ as shown in Fig.\ref{fig:foliationNEW} with a multiplicity of $2$. The cardinality of this volume cover is $\leq \epsilon^{-\delta}$. Next, we construct a \emph{Whitney cover} $\mathcal{W}$ (cf.\ \cite{antti2013,stein1970}) of $K$ by cubes whose diameter is comparable to their distance to $\partial K$ and we restrict this to $\mathcal{W}_\epsilon$ which is initially a collection of elements in $\mathcal{W}$ which have non-zero intersection with
  $K\setminus K_{-\epsilon}$ and adjust this cover s.\,t.\ $\mathcal{W}_\epsilon$ and $\{D_i ^\epsilon\}$ are disjoint. With a variant of the Dirichlet-Neumann bracketing (see Prop.~\ref{prop:bracketing}), we have
  \begin{align*}
      N_N(K,\lambda) \leq \sum_i N_N(D_i ^\epsilon,2\lambda)+ N_N \bigg(\inn \bigcup_{Q \in \mathcal{W}_\epsilon} \overline{Q},\lambda \bigg).
  \end{align*} 
  Now \cite{DOKUMENT} provides an estimate for the first non-trivial Neumann eigenvalue of $D_i ^\epsilon$ showing that it is proportional (denoted by $\sim$) to $\epsilon^{-2}$ so that $N_N(D_i ^\epsilon,\lambda) = 1$ for all $\lambda \leq \lambda_0$ where $\lambda_0 \sim \epsilon^{-2}$. Then for $\epsilon$ sufficiently small, $\sum_i N_N(D_i ^\epsilon,\lambda) = \#\{D_i ^\epsilon\} \leq \epsilon^{-\delta} \sim \lambda^{\delta/2}$.
  Since $\inn \bigcup_{Q \in \mathcal{W}_\epsilon} \overline{Q}$ is a planar polygonal region, it is sufficient to obtain an estimate on its circumference, which is directly related to the number of Whitney cubes in $\mathcal{W}_\epsilon$ of smallest diameter. This circumference turns out to be proportional to $\epsilon^{1-\delta}$. This amounts to a second term proportional to $\epsilon^{1-\delta}\lambda^{1/2} \sim \lambda ^{\delta/2}$.
   \begin{figure}[ht]
 \centering
  \begin{tikzpicture}[
 declare function = {h(\g) = -sqrt(3)*\g + sqrt(3)/3; h2(\var) = -sqrt(3)/(9*6)*\var + 36*sqrt(3)/(9*6); hx(\var) = -sqrt(3)*\var + sqrt(3)/6; },
			scale=5,decoration=Koch snowflake]
\clip (0.6,-0.35) rectangle (2.6,0.9);
    \draw[line width=0.01mm] decorate{decorate{decorate{ decorate{ decorate{ decorate{ (0,0) -- (3,0) }}}}}};
    \draw[fill=black!20!white!40!,line width=0.01mm,dashed] decorate{ (0,0) -- (3,0) };
    \draw[white] (0,0) -- (3,0);
    \draw[fill=black!80!white!40!] decorate{decorate{decorate{ decorate{(1,0) -- (3/2,0.866025) }}}};
    \draw[fill=black!80!white!40!] decorate{decorate{decorate{ decorate{(3-3/2,0.866025) -- (3-1,0) }}}};
    \draw[fill=black!40!white!40!,line width=0.01mm,dashed] decorate{ decorate{(1,0) -- (3/2,0.866025) }};
    \draw[fill=black!40!white!40!,line width=0.01mm,dashed] decorate{ decorate{(3-3/2,0.866025) -- (3-1,0) }};
    \draw[fill=black!40!white!40!,line width=0.01mm] decorate{decorate{ (1,0) -- (1+1/18,0.866025/9) }};
    \draw[fill=black!40!white!40!,line width=0.01mm] decorate{decorate{ (3-1-1/18,0.866025/9) -- (3-1,0) }};
    \draw[fill=black!40!white!40!,line width=0.01mm] decorate{decorate{ (1+0.111111,0.19245) -- (1+1/18+0.111111,0.866025/9+0.19245) }};
    \draw[fill=black!40!white!40!,line width=0.01mm] decorate{decorate{ (3-1-1/18-0.111111,0.866025/9+0.19245) -- (3-1-0.111111,0.19245) }};
    \draw[fill=black!40!white!40!,line width=0.01mm] decorate{decorate{ (1+0.333333, 0.57735) -- (1+1/18+0.333333,0.866025/9+0.57735) }};
    \draw[fill=black!40!white!40!,line width=0.01mm] decorate{decorate{ (3-1-1/18-0.333333,0.866025/9+0.57735) -- (3-1-0.333333, 0.57735) }};
    \draw[fill=black!40!white!40!,line width=0.01mm] decorate{decorate{ (1+0.444444, 0.7698) -- (1+1/18+0.444444,0.866025/9+0.7698) }};
    \draw[fill=black!40!white!40!,line width=0.01mm] decorate{decorate{ (3-1-1/18-0.444444,0.866025/9+0.7698) -- (3-1-0.444444, 0.7698)}};
    \draw[fill=black!80!white!40!,line width=0.01mm] decorate{decorate{ decorate{ (1.16667, 0.288675) -- (1, 0.57735) }}};
    \draw[fill=black!80!white!40!,line width=0.01mm] decorate{decorate{ decorate{ (3-1, 0.57735) -- (3-1.16667, 0.288675)}}};
    \draw[fill=black!80!white!40!,line width=0.01mm] decorate{decorate{ decorate{ (1, 0.57735) -- (1.33333, 0.57735) }}};
    \draw[fill=black!80!white!40!,line width=0.01mm] decorate{decorate{ decorate{ (3-1.33333, 0.57735) -- (3-1, 0.57735) }}};
    \draw[fill=black!120!white!40!,line width=0.01mm] decorate{decorate{ decorate{ (1.11111, 0.3849) -- (1, 0.3849) }}};
    \draw[fill=black!120!white!40!,line width=0.01mm] decorate{decorate{ decorate{ (3-1, 0.3849) -- (3-1.11111, 0.3849)}}};
    \draw[fill=black!120!white!40!,line width=0.01mm] decorate{decorate{ decorate{ (1, 0.3849) -- (1.05556, 0.481125) }}};
    \draw[fill=black!120!white!40!,line width=0.01mm] decorate{decorate{ decorate{ (3-1.05556, 0.481125) -- (3-1, 0.3849) }}};
    \draw[fill=black!120!white!40!,line width=0.01mm] decorate{decorate{ decorate{ (1+0.0555556+0.0555556, 0.3849+0.096225+0.096225) -- (1.05556+0.0555556+0.0555556, 0.481125+0.096225+0.096225) }}};
    \draw[fill=black!120!white!40!,line width=0.01mm] decorate{decorate{ decorate{ (3-1.05556-0.0555556-0.0555556, 0.481125+0.096225+0.096225) -- (3-1-0.0555556-0.0555556, 0.3849+0.096225+0.096225) }}};
    \draw[fill=black!120!white!40!,line width=0.01mm] decorate{decorate{ decorate{ (1.05556+0.0555556+0.0555556, 0.481125+0.096225+0.096225) -- (1.05556+0.0555556+0.0555556+0.0555556, 0.481125+0.096225+0.096225-0.096225)}}};
    \draw[fill=black!120!white!40!,line width=0.01mm] decorate{decorate{ decorate{ (3-1.05556-0.0555556-0.0555556-0.0555556, 0.481125+0.096225+0.096225-0.096225) -- (3-1.05556-0.0555556-0.0555556, 0.481125+0.096225+0.096225)}}};
	\draw[fill=black!10!white!40!, line width=0.05mm] (1,0) -- (2,0) -- (2,-0.2) -- (1,-0.2) -- (1,0);
	\draw[line width=0.4mm] (1,-0.2) -- (2,-0.2);
	\node[scale=1] at (1.4,0.065) {$D_i ^\epsilon$};
    \draw[<->] (2.025,0) -- (2.025,-0.2);
    \node[rotate=90,below] at (2.025,-0.1) {$r=\epsilon$};
    \draw[<->] (1,-0.225) -- (2,-0.225);
    \node[below] at (1.5,-0.225) {$3^{-k}$};
	\foreach \k in {0,1,2,4,5,6,7,8,9,10,11,12,13,14,15,16,17,18,19,20,21,22,23,24,25,26,27,28,29,30,31,33,34,35,36,37,38,39,40,41,42,43,44,45,46,47,48,49,50,51,52,53,54}{.
        \draw (1.5+\k/108,-0.2) -- (1.5+\k/108,{-sqrt(3)*(1.5+\k/108)+2*sqrt(3)});
        }
        \foreach \k in {32}{
        \draw (1.5+\k/108,-0.2) -- (1.5+\k/108,{-sqrt(3)*(1.5+\k/108)+2*sqrt(3)});
        }
        \foreach \k in {1}{
        \draw (1.5+\k/108,{-sqrt(3)*(1.5+\k/108)+2*sqrt(3)}) -- ({1.5+\k/108+sqrt(3)/2*(1/27)*(h2(27+1-\k))},{-sqrt(3)*(1.5+\k/108)+2*sqrt(3)+1/2*(1/27)*(h2(27-1+\k))});
        }
        \foreach \k in {3}{
        \draw (1.5+\k/108,{-sqrt(3)*(1.5+\k/108)+2*sqrt(3)}) -- ({1.5+\k/108+sqrt(3)/2*(1/3)*(1/3)*(h2(27+3-\k))},{-sqrt(3)*(1.5+\k/108)+2*sqrt(3)+1/2*(1/3)*(1/3)*(h2(27-3+\k))});
        }
        \foreach \k in {5}{
        \draw (1.5+\k/108,{-sqrt(3)*(1.5+\k/108)+2*sqrt(3)}) -- ({1.5+\k/108+sqrt(3)/2*(1/27)*(h2(27+5-\k))},{-sqrt(3)*(1.5+\k/108)+2*sqrt(3)+1/2*(1/27)*(h2(27-5+\k))});
        }
        \foreach \k in {6,7,8,9}{
        \draw (1.5+\k/108,{-sqrt(3)*(1.5+\k/108)+2*sqrt(3)}) -- ({1.5+\k/108+sqrt(3)/2*(1/3)*(h2(54-3*\k))},{-sqrt(3)*(1.5+\k/108)+2*sqrt(3)+1/2*(1/3)*(h2(54-3*\k))});
        }
        \foreach \k in {9,10,11,12}{
        \draw (1.5+\k/108,{-sqrt(3)*(1.5+\k/108)+2*sqrt(3)}) -- ({1.5+\k/108+sqrt(3)/2*(1/3)*(h2(3*\k))},{-sqrt(3)*(1.5+\k/108)+2*sqrt(3)+1/2*(1/3)*(h2(3*\k))});
        }
        \foreach \k in {13}{
        \draw (1.5+\k/108,{-sqrt(3)*(1.5+\k/108)+2*sqrt(3)}) -- ({1.5+\k/108+sqrt(3)/2*(1/27)*(h2(27))},{-sqrt(3)*(1.5+\k/108)+2*sqrt(3)+1/2*(1/27)*(h2(27))});
        }
        \foreach \k in {15}{
        \draw (1.5+\k/108,{-sqrt(3)*(1.5+\k/108)+2*sqrt(3)}) -- ({1.5+\k/108+sqrt(3)/2*(1/3)*(1/3)*(h2(27+15-\k))},{-sqrt(3)*(1.5+\k/108)+2*sqrt(3)+1/2*(1/3)*(1/3)*(h2(27-15+\k))});
        }
        \foreach \k in {17}{
        \draw (1.5+\k/108,{-sqrt(3)*(1.5+\k/108)+2*sqrt(3)}) -- ({1.5+\k/108+sqrt(3)/2*(1/27)*(h2(27))},{-sqrt(3)*(1.5+\k/108)+2*sqrt(3)+1/2*(1/27)*(h2(27))});
        }
        \foreach \k in {18,19,20,21,22,23,24,25,26,27}{
        \draw (1.5+\k/108,{-sqrt(3)*(1.5+\k/108)+2*sqrt(3)}) -- ({1.5+\k/108+sqrt(3)/2*(h2(54-\k))},{-sqrt(3)*(1.5+\k/108)+2*sqrt(3)+1/2*(h2(54-\k))});
        }
        \foreach \k in {22,23}{
        \draw ({1.5+\k/108+sqrt(3)/2*(h2(54-\k))},{-sqrt(3)*(1.5+\k/108)+2*sqrt(3)+1/2*(h2(54-\k))}) -- ({1.5+\k/108+sqrt(3)/2*(h2(54-\k))+0},{-sqrt(3)*(1.5+\k/108)+2*sqrt(3)+1/2*(h2(54-\k))+1/9*sqrt(3)/3});
        }
	\foreach \k in {27,28,29,30,31,33,34,35,36}{
        \draw (1.5+\k/108,{-sqrt(3)*(1.5+\k/108)+2*sqrt(3)}) -- ({1.5+\k/108+sqrt(3)/2*(h2(\k))},{-sqrt(3)*(1.5+\k/108)+2*sqrt(3)+1/2*(h2(\k))});
        }
        \foreach \k in {32}{
        \draw (1.5+\k/108,{-sqrt(3)*(1.5+\k/108)+2*sqrt(3)}) -- ({1.5+\k/108+sqrt(3)/2*(h2(\k))},{-sqrt(3)*(1.5+\k/108)+2*sqrt(3)+1/2*(h2(\k))});
        }
        \foreach \k in {31}{
        \draw ({1.5+\k/108+sqrt(3)/2*(h2(\k))},{-sqrt(3)*(1.5+\k/108)+2*sqrt(3)+1/2*(h2(\k))}) -- ({1.5+\k/108+sqrt(3)/2*(h2(\k))+sqrt(3)/2*1/9*sqrt(3)/3},{-sqrt(3)*(1.5+\k/108)+2*sqrt(3)+1/2*(h2(\k))-1/2*1/9*sqrt(3)/3});
        }
        \foreach \k in {32}{
        \draw ({1.5+\k/108+sqrt(3)/2*(h2(\k))},{-sqrt(3)*(1.5+\k/108)+2*sqrt(3)+1/2*(h2(\k))}) -- ({1.5+\k/108+sqrt(3)/2*(h2(\k))+sqrt(3)/2*1/9*sqrt(3)/3},{-sqrt(3)*(1.5+\k/108)+2*sqrt(3)+1/2*(h2(\k))-1/2*1/9*sqrt(3)/3});
        }
        \foreach \k in {37}{
        \draw (1.5+\k/108,{-sqrt(3)*(1.5+\k/108)+2*sqrt(3)}) -- ({1.5+\k/108+sqrt(3)/2*(1/27)*(h2(27))},{-sqrt(3)*(1.5+\k/108)+2*sqrt(3)+1/2*(1/27)*(h2(27))});
        }
        \foreach \k in {39}{
        \draw (1.5+\k/108,{-sqrt(3)*(1.5+\k/108)+2*sqrt(3)}) -- ({1.5+\k/108+sqrt(3)/2*(1/3)*(1/3)*(h2(27+39-\k))},{-sqrt(3)*(1.5+\k/108)+2*sqrt(3)+1/2*(1/3)*(1/3)*(h2(27-39+\k))});
        }
        \foreach \k in {41}{
        \draw (1.5+\k/108,{-sqrt(3)*(1.5+\k/108)+2*sqrt(3)}) -- ({1.5+\k/108+sqrt(3)/2*(1/27)*(h2(27))},{-sqrt(3)*(1.5+\k/108)+2*sqrt(3)+1/2*(1/27)*(h2(27))});
        }
        \foreach \k in {42,43,44,45}{
        \draw (1.5+\k/108,{-sqrt(3)*(1.5+\k/108)+2*sqrt(3)}) -- ({1.5+\k/108+sqrt(3)/2*(1/3)*(h2(162-3*\k))},{-sqrt(3)*(1.5+\k/108)+2*sqrt(3)+1/2*(1/3)*(h2(162-3*\k))});
        }
        \foreach \k in {45,46,47,48}{
        \draw (1.5+\k/108,{-sqrt(3)*(1.5+\k/108)+2*sqrt(3)}) -- ({1.5+\k/108+sqrt(3)/2*(1/3)*(h2(3*\k-108))},{-sqrt(3)*(1.5+\k/108)+2*sqrt(3)+1/2*(1/3)*(h2(3*\k-108))});
        }
        \foreach \k in {49}{
        \draw (1.5+\k/108,{-sqrt(3)*(1.5+\k/108)+2*sqrt(3)}) -- ({1.5+\k/108+sqrt(3)/2*(1/27)*(h2(27))},{-sqrt(3)*(1.5+\k/108)+2*sqrt(3)+1/2*(1/27)*(h2(27))});
        }
        \foreach \k in {51}{
        \draw (1.5+\k/108,{-sqrt(3)*(1.5+\k/108)+2*sqrt(3)}) -- ({1.5+\k/108+sqrt(3)/2*(1/3)*(1/3)*(h2(27+51-\k))},{-sqrt(3)*(1.5+\k/108)+2*sqrt(3)+1/2*(1/3)*(1/3)*(h2(27-51+\k))});
        }
        \foreach \k in {53}{
        \draw (1.5+\k/108,{-sqrt(3)*(1.5+\k/108)+2*sqrt(3)}) -- ({1.5+\k/108+sqrt(3)/2*(1/27)*(h2(27))},{-sqrt(3)*(1.5+\k/108)+2*sqrt(3)+1/2*(1/27)*(h2(27))});
        }
	\foreach \k in {3}{
        \draw (1.5+\k/108,-0.2) -- (1.5+\k/108,{-sqrt(3)*(1.5+\k/108)+2*sqrt(3)});
        }
\end{tikzpicture}
\caption{Instance of a volume covering domain as used in the proof of Thm.~\ref{thm:satz1}.}\label{fig:foliationNEW}
\end{figure}
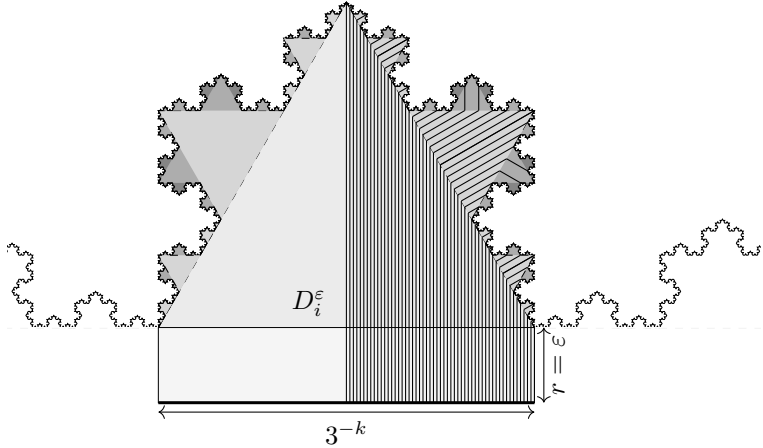
  \end{proof}
  \begin{rem}
    The general strategy behind the proof of Thm.~\ref{thm:satz1} is similar to the concept shown in \cite{NetrusovSafarov2005}, however we are not limited to domains whose domains are locally graphs.\par
    In \cite{la1991}, Lapidus found a similar asymptotic result under the condition that the upper Minkowski content exists.
    In contrast, the result of Thm.~\ref{thm:satz1} is based on an estimate on the Neumann counting function which only relies on the existence of the upper inner Minkowski content and provides explicit upper bounds. Moreover, the explicit bounds in \cite{DOKUMENT} are given for all $\lambda\geq 0.1$ rather than asymptotically. This is important for the application below.
\end{rem}
  \section{Asymptotics of counting functions}\label{sec:weylasymp}
  Let $K \subset \mathbb{R}^2$ be the Koch snowflake with base length $1$ and let $\Omega(k_1,k_2)$ be the limit set described in Fig.\subref{fig:kochflockesec5}-\subref{fig:kochflockesec5k1k2} and Sec.\ref{sec:kochspray} with corresponding generator $K$.
    \begin{thm}
      There is an asymptotic expansion of $N_D(\Omega(k_1,k_2),\textup{e}^t)$ of the form
      \begin{align*}
          N_D(\Omega(k_1,k_2),{\textup{e}}^t) = \frac{1}{4\pi} \vol_2(\Omega(k_1,k_2)) {\textup{e}}^{t} - \sum_{\substack{z \in \mathcal{Z}_{\textup{C}},\\ \Re(z)<-\delta/2}} \widetilde{{\coeffweyl}_{\beta(t)}}(z) {\textup{e}}^{-2a\lfloor \frac{t}{2a} \rfloor z} + {o}({\textup{e}}^{t(\delta/2+\gamma)})
      \end{align*}
      as $t \to \infty$ for any $\gamma>0$. The absolute values of $\widetilde{{\coeffweyl}_{\beta(t)}}(z)$ are bounded from above according to Tab.~\ref{tab:upperboundsG}. Here, $\beta(t) := 2a\{ \frac{t}{2a}\}$, with $\{x\}$ denoting the fractional part of $x$, and  
      $\mathcal{Z}_{\text{C}} :=\big{\{} z\in\mathbb C \ :\  \sum_{i \in \Sigma} r_i^{-2z} =1,\ \Im(z)\in\big{[}0,\frac{\pi}{a}\big{)} \big{\}}$.
  \end{thm}
  \begin{proof} We first consider $\Omega:= \Omega(0,0)$. We define $N(t) := N_D(\Omega,{\textup{e}}^t)$ and $g(t):=N_D(K,{\textup{e}}^t) = \frac{1}{4\pi} \vol_2(K) {\textup{e}}^{t} +M(t){\textup{e}}^{t\delta/2}$ for an $M \in  \mathcal{O}(1)$ as shown in Thm.\ref{thm:satz1}. By \ref{item:2},\ref{item:3} in Sec.\ref{sec:preliminaries} and \eqref{eq:Omega}, $N(t)$ satisfies
  \begin{align*}
      N(t) &= \sum_{k \geq 0} \sum_{w \in \Sigma^k} N_D(\phi_w K,\textup{e}^t)\\
           &= \sum_{k \geq 0 } \sum_{w \in \Sigma^k} g(t - 2\log r_{w_1}-2\log r_{w_2}  \cdots -2\log r_{w_k}) .
  \end{align*}
  Let $\ell_0 \in \mathbb{Z}$ be maximal such that $N(2a\ell_0)=0$. Then for every $\beta \in [0,2a)$ we define its two-sided Fourier-Laplace transform $\widehat{N}_\beta(z)$ for $\Re(z) < -1$ and rewrite this to isolate the poles of its maximal meromorphic extension. We make use of the fact that $\{-2\log r_{i}\,:\,i \in \Sigma\} \subset 2a\mathbb{Z}$ which follows from \eqref{eq:nu}, so that for any $w \in \Sigma^k$ there is an $\nu_{\omega} \in \mathbb{Z}$ with $-2\log r_{\omega} = 2a \nu_{\omega}$, and perform an index shift $\ell \to \ell + \sum_{i=1} ^k \nu_i = \widetilde{\ell}$.
  \begin{align*}
      \widehat{N}_\beta(z) &:= \sum_{\ell \in \mathbb{Z}} {\textup{e}}^{2a\ell z}N(2a\ell+\beta) = \sum_{k \geq 0 } \sum_{w \in \Sigma^k} \sum_{\ell \in  \mathbb{Z}} {\textup{e}}^{2a\ell z} g(2a\ell + \beta - 2\log r_{w_1} \cdots -2\log r_{w_k} )\\
      &= \sum_{k \geq 0 } \sum_{w \in \Sigma^k}{\textup{e}}^{z\left( - 2\log r_{w_1}- \cdots -2\log r_{w_k} \right)} \sum_{\widetilde{\ell} \in  \mathbb{Z}} {\textup{e}}^{2a\widetilde{\ell} z} g(2a\widetilde{\ell} + \beta)
      \\
      &=\sum_{k \geq 0 } \sum_{w \in \Sigma^k} {r_w}^{-2z} \left( \sum_{\ell \geq \ell_0} {\textup{e}}^{2a\ell z} \frac{1}{4\pi} \vol_2(K) {\textup{e}}^{2a\ell+\beta} + \sum_{\ell \geq \ell_0} {\textup{e}}^{2a\ell z} M(2a\ell+\beta){\textup{e}}^{(2a\ell+\beta)\delta/2} \right)\\
      &= \underbrace{\frac{1}{1-\sum_{i \in \Sigma}{r_i}^{-2z}}}_{:= P(z)} \underbrace{\left(\frac{\frac{1}{4\pi} \vol_2(K){\textup{e}}^{\beta} {\textup{e}}^{2a\ell_0(z+1)}}{1-{\textup{e}}^{2a(z+1)}} + M_\beta(z) \right)}_{:={\coeffweyl}_\beta(z)}
  \end{align*}
  for some complex function $M_\beta$ bounded in $\mathbb{C}_{-\delta / 2} := \{z \in \mathbb{C}:\Re(z) < -\delta/2\}$ by
  \begin{align*}
      0 \leq |M_\beta(z)| \leq \left| \sum_{\ell \geq \ell_0} {\textup{e}}^{2a\ell z} M(2a\ell+\beta){\textup{e}}^{(2a\ell+\beta)\delta/2} \right| \leq \left| \frac{ \widetilde{M_K} {\textup{e}}^{\beta \delta/2} }{1-{\textup{e}}^{2a(z+\delta/2)}} \right|,
  \end{align*}
  where $\widetilde{M_K} := \max(C_-,C_+)$ is taken from Thm.~\ref{thm:satz1}. 
  Therefore, $\widehat{N}_{\beta}$ can be meromorphically extended to $\mathbb C_{-\delta/2}$.
  Let 
  \begin{align*}
    \mathcal{Z}_{\text{C}} :=\bigg{\{} z\in\mathbb C \ :\  \sum_{i \in \Sigma} r_i^{-2z} =1,\ \Im(z)\in\big{[}0,\frac{\pi}{a}\big{)} \bigg{\}}
    \end{align*}
    Then, since all poles in $\mathcal{Z}_C$ are simple and $-1 \notin \mathcal{Z}_{\textup{C}}$,
  \begin{align*}
      &\sum_{\ell \geq 0} {\textup{e}}^{2a\ell s} N(2a\ell+\beta) - \sum_{z \in \mathcal{Z}_{\textup{C}}\cap\mathbb C_{-\delta/2}} \frac{-2a{\coeffweyl}_\beta(z) \Res_z P}{1 - {\textup{e}}^{2a(s-z)}} - \left( \frac{\frac{1}{4\pi} \vol_2(K) {\textup{e}}^{\beta} }{1-\sum_{i \in \Sigma} {r_i}^{2}} \right) \frac{1}{1-{\textup{e}}^{2a(s+1)}}
      \\
      &= \sum_{\ell \geq 0} {\textup{e}}^{2a\ell s} \left( N(2a\ell+\beta) + \sum_{z \in \mathcal{Z}_{\textup{C}}\cap\mathbb C_{-\delta/2}} \frac{-2a{\coeffweyl}_\beta(z)}{\sum_{i\in \Sigma} 2\log (r_i){r_i}^{-2z}}{\textup{e}}^{-2a\ell z} - \frac{1}{4\pi}\cdot \frac{\vol_2(K) }{1-\sum_{i \in \Sigma} {{r_i}}^{2} } {\textup{e}}^{2a\ell + \beta}  
      \right)\\
      &=: \sum_{\ell \geq 0} {\textup{e}}^{2a\ell s} \cdot Y_{\ell}
  \end{align*}
  is holomorphic in $\mathbb{C}_{-\delta / 2}$. 
  As a power series has a singularity on its radius of convergence, and the above series expansion is holomorphic in $\mathbb{C}_{-\delta / 2}$ we deduce that $ Y_{\ell} \in {o}({\textup{e}}^{2a\ell(\delta/2+\gamma)})$ for any $\gamma>0$ as $\ell \to \infty$. Thus, for any $\gamma>0$  
\begin{align*}
    N(2a\ell+\beta) 
    &= 
    \frac{1}{4\pi} \cdot \frac{\vol_2(K) }{1-\sum_{i \in \Sigma} {{r_i}}^{2} } {\textup{e}}^{2a\ell + \beta}
    - \sum_{z \in \mathcal{Z}_{\textup{C}}\cap\mathbb C_{-\delta/2}} \frac{-2a {\coeffweyl}_\beta(z)}{\sum_{i\in \Sigma} 2\log (r_i){r_i}^{-2z}}{\textup{e}}^{-2a\ell z} 
    + {o}({\textup{e}}^{a\ell(\delta+2\gamma)})
\end{align*}
as $\ell \to \infty$.
  Notice that $\frac{\vol_2(K) }{1-\sum_{i \in \Sigma} {{r_i}}^{2} } = \vol_2(K) \sum_{\omega \in \Sigma^\ast} {r_{\omega}}^{2} = \vol_2(\Omega)$ as expected. Writing $t = 2a\lfloor \frac{t}{2a} \rfloor + \beta(t)$ with $\beta(t)=2a\{\frac{t}{2a}\}$, this implies an asymptotic expansion of $N(t)$ given by
  \begin{align}
      N(t) = \frac{1}{4\pi} \vol_2(\Omega) {\textup{e}}^{t} - \sum_{z \in \mathcal{Z}_{\textup{C}}\cap\mathbb C_{-\delta/2}} \widetilde{{\coeffweyl}_{\beta(t)}}(z) {\textup{e}}^{-2a\lfloor \frac{t}{2a} \rfloor z} + {o}({\textup{e}}^{t(\delta/2+\gamma)}), \label{eq:Ntgl}
  \end{align}
  with $\widetilde{{\coeffweyl}_{\beta(t)}}(z):= \frac{-2a{\coeffweyl}_{\beta(t)}(z)}{\sum_{i\in \Sigma} 2\log ({r_i}){r_i}^{-2z}}$. Since $t \mapsto {\coeffweyl}_{\beta(t)}(z)$ is bounded, the asymptotic expansion in \eqref{eq:Ntgl} already contains all terms corresponding to a growth rate ${\textup{e}}^{tz'}$ with $z' \in (\delta/2,1)$.\par
  This estimation remains correct in the case of $(k_1,k_2) \neq (0,0)$ after adapting $C_\pm$ from Thm.~\ref{thm:satz1}. By \ref{item:3} in Sec.~\ref{sec:preliminaries}, $N_D(\alpha \Omega,t) = N_D(\Omega,\alpha^2 t) \leq C_W \vol_2(\alpha \Omega)t + C_+ \alpha^{\delta}t^{\delta/2} $. So we obtain an expression for an upper bound of the remainder term corresponding to $\Omega(k_1,k_2)$ if we replace $C_\pm$ with $C_\pm + 9^{-\delta} k_1 C_\pm + (9\sqrt{3})^{-\delta} k_1 C_\pm$.
 In Tab.~\ref{tab:upperboundsG}, we provide approximates of the exponents and bounds on the coefficients corresponding to $z \in \mathcal{Z}_{\text{C}} \cap (-1,-\delta/2)$ for three different allocations of $(k_1,k_2)\in\{0,\ldots,6\}^2$ exhibiting different arrangements of the relevant poles.
\begin{table}[h]
 \centering
 \begin{tabular}{c|l|l}
 $(k_1,k_2)$ & Approx. Values of $z \in \mathcal{Z}_{\textup{C}}$ with $\Re (z) < -\frac{\delta}{2}$ up to $\frac{2\pi \mathbf{i}}{a}\mathbb{Z}$ & Upper bound of $|\widetilde{{\coeffweyl}_{\beta(t)}}(z)|$\\
 \hline
 $(0,0)$ & $-0.952455$ & $1.68\cdot 10^6$\\
 $(0,6)$ & $-0.928326$, $-0.71134 \pm 2.58082\mathbf{i}$ & $1.81 \cdot 10^6$, $2.45 \cdot 10^5$\\
 $(6,6)$ & $-0.888243$, $-0.839089 \pm 1.34671\mathbf{i}$, $-0.666227 \pm 2.8596\mathbf{i}$ & $1.68\cdot 10^6$, $3.46\cdot 10^5$, $2.92\cdot 10^5$
\end{tabular}
\caption{Upper bounds for coefficients of the asymptotic expansion at the relevant poles.}\label{tab:upperboundsG}
\end{table}
\end{proof}
\begin{rem}
      A critical problem occurs in the asymptotic expansion in the following sense: Suppose $N_D(\Omega,t) = { (2\pi)^n V_n } t^{n/2} + M(\log t)t^{\delta/2}+A(\log t)t^{\delta/2}$ with some bounded $M$ and
  \begin{align*}
      A: t \mapsto 
      \begin{cases}
          \frac{1}{\lfloor t/a\rfloor} &\text{ if }\exists m \in \mathbb{N}_0:\lfloor t/a\rfloor = 2^m\\
          0 &\text{ else }
      \end{cases}
  \end{align*}
  so that $A(\log t) \in o(1)$ and $g(t)={ (2\pi)^n V_n } {\textup{e}}^{tn/2} + M(t){\textup{e}}^{t\delta/2}+A(t){\textup{e}}^{t\delta/2}$. Then this third term in $g(t)$ leads to the following term in the two-sided Fourier-Laplace transform $\widehat{g}_\beta(z)$:
  \begin{align*}
      \sum_{m \in \mathbb{N}_0} \frac{\left( {\textup{e}}^{a(z+\delta/2)} \right)^{2^m}}{2^m}.
  \end{align*}
  If this had a meromorphic extension beyond $\Re(z)<-\delta/2$, so would $\widehat{g}_\beta '(z)$. But $\widehat{g}_\beta '(z)$ contains a series of the form
  \begin{align*}
      \sum_{m \in \mathbb{N}_0} \left( {\textup{e}}^{a(z+\delta/2)} \right)^{2^m},
  \end{align*}
  which diverges whenever $\e^{a(z+\delta/2)}$ is a root of unity of any power of $2$. In other words, whenever $z = -\frac{\delta}{2} + \frac{2\pi \mathbf{i} q}{a 2^p}$ for some $p \in \mathbb{Z}$ and $q \in \mathbb{N}$. However, this set lies dense in $-\frac{-\delta}{2}+\mathbf{i}\mathbb{R} \subset \mathbb{C}$. This shows that further information about the behaviour of the remainder-term of $N_D(\Omega,t)$ is necessary in order to ensure existence of a meromorphic extension of $\widehat{g}_\beta(z)$ that has a pole at $-\delta/2$.
  \end{rem}

   \section{Asymptotics of parallel volumes}\label{sec:minkcontent}
   In this section, we will derive an asymptotic expansion of $\e\mapsto\vol_2\left(\Omega(k_1,k_2)_{-\e}\right)$ for the sets $\Omega(k_1,k_2)$ with $(k_1,k_2)\in\{0,\ldots,6\}^2$ that we defined in Sec.~\ref{sec:kochspray}.
   \begin{thm}\label{thm:volume}
       Let $\Omega(k_1,k_2)$ be as in Sec.~\ref{sec:kochspray}. Then for all $\beta\in[0,a)$ we have an asymptotic expansion of $\vol_2 \left( \Omega(k_1,k_2)_{-\textup{e}^{-(a\ell+\beta)}}\right)$ as $\ell\to\infty$ of the following form. 
       \begin{align*}
           \vol_2 \left( \Omega(k_1,k_2)_{-\textup{e}^{-(a\ell+\beta)}}\right)
     = {\coeffmink}_{\beta}(2)\textup{e}^{-2a\ell} 
     + {\coeffmink}_{\beta}\left(2- \frac{\log 4}{\log 3}\right)\textup{e}^{-a\ell(2-\frac{\log 4}{\log 3})} 
     +\sum_{z\in \mathcal{Z}_{\textup P}} \textup{e}^{-a\ell z}{\coeffmink}_{\beta}(z)
     +o(\textup{e}^{-a\ell\gamma}),
       \end{align*}
       as $\ell\to\infty$ for any $\gamma>0$,
       with coefficients ${\coeffmink}_{\beta}$ given in \eqref{eq:Minkcoeff} and evaluated in Tab.~\ref{tab:valuecontentG}. Here, $\mathcal{Z}_{\text{P}} := \{ z\in\mathbb C \ :\  \sum_{i \in \Sigma} r_i^{2-z} =1,\ \Im(z)\in\big{[}0,\frac{2\pi}{a}\big{)} \}$.
   \end{thm}
   A key part of the proof of Thm.~\ref{thm:volume} relies on precise knowledge of the inner $\e$-parallel volume of the generator $\gen:= O\setminus \overline{\Phi(k_1,k_2) O}$ for all $\e>0$. As $\gen$ is a disjoint union of Koch snowflakes of different sizes, a key step in the proof is to determine the inner $\e$-parallel volume of the Koch snowflake of base length 1 in the next lemma.
   
\begin{lem}\label{lem:innerKochflake}
Let $K$ denote the filled-in Koch snowflake with base-length $1$. The map $\e\mapsto\vol_2(K_{-\e})$, defined on the positive reals, which maps $\epsilon$ to the area of the inner $\e$ neighbourhood $K_{-\epsilon}:=\{x\in K\ :\ \inf_{y\in \partial K}\|x-y\|_2\leq \epsilon\}$, is continuous and can be evaluated as follows.
\begin{align}\label{eq:innerkochvolume}
    \vol_2(K_{-\e}) =
    \begin{cases}
        \frac{2\sqrt{3}}{5} & :\e>\frac{1}{3}\\
        \frac{7\sqrt{3}}{30} + \sqrt{\e^2-\frac{1}{36}} + 6\e^2\arcsin{\left(\frac{1}{6\e}\right)}-\pi\e^2 & :\frac{\sqrt{3}}{9} <\e \leq \frac{1}{3}\\
        \frac{8\sqrt{3}}{45} + \pi\e^2 + 12\cdot\vol_2( K_{-\e}\cap \Gamma)& :\frac{1}{9}<\e\leq\frac{\sqrt{3}}{9}\\
        u\circ\alpha(\e)\e^{2-\log{4}/\log{3}} + v\circ\alpha(\e)\e^2 & : \e\leq\frac{1}{9}\ \text{and}\ \alpha(\e)<\frac{1}{2}\\
        \widetilde{u}\circ\alpha(\e)\e^{2-\log{4}/\log{3}} + v\circ\alpha(\e)\e^2 & : \e\leq\frac{1}{9}\ \text{and}\ \alpha(\e)\geq\frac{1}{2}.
    \end{cases}
\end{align}
Here, $\{t\}:= t - \lfloor t\rfloor$ denotes the fractional part of a real number $t$, $\alpha(\e):=\left\{ -\frac{\log{\e}}{\log{3}}\right\}$, $\Gamma$ is the equilateral filled-in triangle shown in Fig.~\subref{fig:Gammatilde} and defined in the proof of this lemma. Further, $u,\widetilde{u},v\colon [0,1)\to \mathbb R$ are given by
\begin{align*}
    u(t)&:= 
    \left(\frac{9}{4}\right)^{t}\cdot\left[
    \frac{21\sqrt{3}}{40}
    +\frac{3}{4}\cdot\sqrt{3^{-2t}-\frac{1}{4}}
    + 81\cdot \vol_2(K_{-3^{-t-2}}\cap \Gamma)
    \right] +\left(\frac14\right)^{t}\cdot\left[ 
    \frac{3}{2}\cdot \arcsin\left(\frac{3^{t}}{2} \right)
    - \frac{\pi}{6}
    \right],    \\
    \widetilde{u}(t)&:= 
    \left(\frac{9}{4}\right)^{t}\cdot\left[
    \frac{2\sqrt{3}}{5}
    + 27\cdot \vol_2(K_{-3^{-t-1}}\cap \Gamma)
    + 81\cdot \vol_2(K_{-3^{-t-2}}\cap \Gamma)
    \right] 
    +\left(\frac14\right)^{t}\cdot \frac{\pi}{3}, \\
    v(t)&:= -\frac{\pi}{3}-324\cdot 9^{t}\cdot \vol_2(K_{-3^{-t-2}}\cap \Gamma).
\end{align*}
   \end{lem}
   \begin{rem}
   Note that for $\e\leq 1/9$ the area of the inner $\e$ neighbourhood $\vol_2(K_{-\e})$ of the filled-in Koch snowflake $K$ has been determined in \cite{LapidusPearse}, where $u\circ\alpha(\e)$, $\widetilde{u}\circ\alpha(\e)$ and $v\circ\alpha(\e)$ are expressed as infinite complex series. With     Lem.~\ref{lem:innerKochflake} we provide a more geometric representation of $\vol_2(K_{-\e})$ and an alternative and simpler proof of its scaling behaviour. 
   \end{rem}
\begin{proof}[Proof of Lem.~\ref{lem:innerKochflake}.]
    Let $F$ be the Koch curve that is generated by the four contractions $\psi_i\colon\mathbb R^2\to\mathbb R^2$ for $i\in\{1,\ldots,4\}$ given by
    \begin{align*}
        \psi_1(x)&=\frac13 x,\ 
        \psi_2(x)=\frac13 %
        R_{\pi/3}(x)+\frac13\begin{pmatrix}
            1 \\ 0
        \end{pmatrix},\
        \psi_3(x)=\frac13 %
        R_{-\pi/3} (x)+\frac16\begin{pmatrix}
            3 \\ \sqrt 3
        \end{pmatrix},\
        \psi_4(x)=\frac13x +\frac13\begin{pmatrix}
            2 \\ 0
        \end{pmatrix},
    \end{align*}
    with $R_\alpha$ denoting the rotation matrix to the angle $\alpha$ about the origin. 
    Further, let $V$ denote the open region bounded by $F$ and the two line segments 
    $\big{\{} t\binom{3}{-\sqrt{3}}\,:\, t\in[0,\frac16] \big{\}}$ and 
    $\big{\{} \binom{1}{0} -t \binom{3}{-\sqrt{3}}\,:\, t\in[0,\frac16] \big{\}}$, see Fig.~\subref{fig:V}. 
    Suppose without loss of generality that the position of $K$ in $\mathbb R^2$ is so that $K=\overline{V\cup V^1 \cup V^2}$ with $V^k$ denoting the image of $V$ under the rotation around 
    $\frac{1}{6}\binom{3}{-\sqrt{3}}$
    by the angle $\frac{2\pi k}{3}$.
    \begin{figure}
    \begin{subfigure}[t]{.48\textwidth}
    \centering
    \begin{tikzpicture}[scale=2.2,decoration=Koch snowflake]
        \draw[line width=0.01mm, fill=gray!20!white] decorate{decorate{decorate{ decorate{ decorate{ decorate{ (0,0) -- (3,0) }}}}}} -- ({1.5},{-sqrt(3)/2}) -- (0,0);
        \node at (1.5,-0.25) {$V$};
        \draw[dashed] (0,0) -- (3,0);
        \draw[<->] (0:0.5) arc(0:-30:0.5) node[midway,right]{{\small{$30^\circ$}}};
    \end{tikzpicture}
    \caption{Visualisation of the set $V$.}
    \label{fig:V}
    \end{subfigure}
    \begin{subfigure}[t]{.48\textwidth}
        \centering
        \begin{tikzpicture}[scale=2.2,decoration=Koch snowflake]
    \draw[line width=0.01mm] decorate{decorate{decorate{ decorate{ decorate{ decorate{ (0,0) -- (3,0) }}}}}};
    \draw (0,0) -- ({1.5},{-sqrt(3)/2}) -- (3,0);
    \draw[fill=white!20!gray] (1,{-sqrt(3)/3}) -- (1,0) -- ({1/2},{-sqrt(3)/6}) ;
    \draw[fill=white!20!gray] (2,0)  -- (2,{-sqrt(3)/3}) -- ({5/2},{-sqrt(3)/6}) -- (2,0);
    \draw[fill=white!20!gray] (1,0) -- ({3/2},{-sqrt(3)/6}) -- ({3/2},{sqrt(3)/6}) -- (1,0);
    \draw[fill=white!20!gray] (2,0) -- ({3/2},{-sqrt(3)/6}) -- ({3/2},{sqrt(3)/6}) -- (2,0);
    \draw[pattern=north east lines, pattern color= gray] (1,0) -- (1,{-sqrt(3)/3}) -- ({1.5},{-sqrt(3)/2}) -- ({3/2},{-sqrt(3)/6});
    \draw[pattern=north west lines, pattern color= gray] (2,0) -- (2,{-sqrt(3)/3}) -- ({1.5},{-sqrt(3)/2}) -- ({3/2},{-sqrt(3)/6});
    \draw (3/2,{sqrt(3)/2}) -- (3/2,{-sqrt(3)/2});
    \node at (0.8,-0.3) {$\Gamma$};
    \node at (1.2,-0.4) {$\Lambda$};
    \node at (0.5,-1/10) {\small$\psi_1V$};
    \node at (1.3,0.45) {\small{$\psi_2V$}};
    \node at (1.7,0.45) {\small{$\psi_3V$}};
    \node at (2.5,-1/10) {\small{$\psi_4V$}};
    \end{tikzpicture}
    \caption{Decomposition of $V$ into four congruent copies of $\Gamma$, two congruent copies of $\Lambda$ and the sets $\psi_1(V),\ldots,\psi_4(V)$.}
    \label{fig:Gammatilde}
    \end{subfigure}
\end{figure}
    Then
    \begin{equation}
        \vol_2\left({K}_{-\e}\right) = 3\cdot \vol_2\left({K}_{-\e}\cap V\right).
    \end{equation}
    Therefore, in the following, we focus on determining $\vol_2({K}_{-\e}\cap V)$. For this, let $\Gamma$ denote the filled-in equilateral triangle with vertices 
    $\frac{1}{18}\binom{3}{-\sqrt{3}}$, 
    $\frac{1}{9}\binom{3}{-\sqrt{3}}$ and 
    $\frac{1}{3}\binom{1}{0}$. 
    Moreover, $\Lambda$ will denote the filled-in rhombus with vertices 
    $\frac{1}{3}\binom{1}{0}$, 
    $\frac{1}{9}\binom{3}{-\sqrt{3}}$, 
    $\frac{1}{6}\binom{3}{-\sqrt{3}}$ and 
    $\frac{1}{18}\binom{9}{-\sqrt{3}}$. 
    The sets $\Gamma$ and $\Lambda$ are depicted in Fig.~\subref{fig:Gammatilde}. For large enough $\e$, the sets $\Lambda$, $\Gamma$ and $\psi_i(V)$, $i\in\{1,\ldots,4\}$, are fully contained in ${K}_{-\e}$. This changes when $\e$ decreases (see Fig.~\subref{fig:cases-proof}) and in the following we distinguish between different cases corresponding to $\e$, where this behaviour changes.

    \emph{Case 1:} If $\e>\frac13$ then
    \begin{align*}
        \vol_2\left({K}_{-\e}\cap V\right) 
        = \vol_2\left( V\right)
        = \vol_2\left( V\cap \mathbb R\times \mathbb R^+\right) + \vol_2\left( V\cap \mathbb R\times \mathbb R^-\right)
        = \frac{\sqrt 3}{20} +\frac{\sqrt 3}{12} 
        = \frac{2\sqrt 3}{15}.
    \end{align*}

    \emph{Case 2:} If $\e\in(\frac{\sqrt 3}{9},\frac13]$ then
    \begin{align*}
        \vol_2\left({K}_{-\e}\cap V\right) 
        &= \sum_{i=1}^4\vol_2\left( {K}_{-\e}\cap \psi_i(V)\right) + 4\vol_2\left( {K}_{-\e}\cap \Gamma\right)
        + 2\vol_2\left( {K}_{-\e}\cap \Lambda\right)\\
        &= 4\vol_2\left( \psi_1(V)\right) + 4\vol_2( \Gamma)+2\vol_2\left( {K}_{-\e}\cap \Lambda\right).
    \end{align*}
    With $\rho=\rho(\e)=\frac{2\pi}{3}-2\arccos\left(\frac1{6\e}\right)$ being the angle that is shown in Fig.~\subref{fig:Case2}, we can evaluate $\vol_2\left( {K}_{-\e}\cap \Lambda\right)$ as follows.
    \begin{align*}
        \vol_2\left( {K}_{-\e}\cap \Lambda\right)
        & = \frac{\rho(\e)}{2}\e^2 + 2\cdot \frac12 \cdot \frac{\sqrt 3}{9}\cdot\e\cdot\sin\left(\frac{\pi}{6}-\frac{\rho(\e)}{2}\right) \\
        & = \frac{\pi}3 \e^2 
        - \e^2\arccos\left( \frac{1}{6\e}\right) 
        - \frac{\sqrt 3}{108}
        +  \frac16\sqrt{\e^2-\frac{1}{36}} 
    \end{align*}
    It follows that 
    \begin{align*}
        \vol_2\left({K}_{-\e}\cap V\right) 
        &= \frac49\cdot\frac{2\sqrt{3}}{15} + 4\cdot \frac{\sqrt{3}}{108}+
        \frac23 \pi\e^2 
        - 2\e^2\arccos\left( \frac{1}{6\e}\right)
        - \frac{\sqrt 3}{54}
        +  \frac13\sqrt{\e^2-\frac{1}{36}}\\
        &= \frac{7\sqrt{3}}{90} + 
        \frac23 \pi\e^2 
        - 2\e^2\arccos\left( \frac{1}{6\e}\right)
        +  \frac13\sqrt{\e^2-\frac{1}{36}}.
    \end{align*}
    Using the identity $\arccos(x)=\frac{\pi}{2}-\arcsin(x)$ the assertion follows.

    \emph{Case 3:} If $\e\in(\frac19,\frac{\sqrt 3}{9}]$ then
    \begin{align*}
        \vol_2\left({K}_{-\e}\cap V\right) 
        &= \sum_{i=1}^4\vol_2\left( {K}_{-\e}\cap \psi_i(V)\right) + 4\cdot\vol_2\left( {K}_{-\e}\cap \Gamma\right)
        + 2\cdot\vol_2\left( {K}_{-\e}\cap \Lambda\right)\\
        &= \frac{4}{9}\cdot \frac{2\sqrt{3}}{15} + 4\vol_2( {K}_{-\e}\cap\Gamma) + \frac{\pi}{3}\e^2.
    \end{align*}
    
    \emph{Case 4:} Suppose that $\e\leq\frac{1}{9}$. Let $W:=V\setminus \bigcup_{i=1}^4 \psi(V)$. Then 
    \begin{equation*}
    V=\bigcup_{k=0}^{\infty} \bigcup_{\omega\in\{1,\ldots,4\}^k} \psi_{\omega}(W) \cup \bigcap_{k=0}^{\infty} \bigcup_{\omega\in\{1,\ldots,4\}^k} \psi_{\omega}(V),
    \end{equation*}
    where all unions are disjoint. As $\bigcap_{k=0}^{\infty} \bigcup_{\omega\in\{1,\ldots,4\}^k} \psi_{\omega}(V) \subset \overline{\bigcap_{k=0}^{\infty} \bigcup_{\omega\in\{1,\ldots,4\}^k} \psi_{\omega}(V)} \subset F$, we know that $\vol_2(\bigcap_{k=0}^{\infty} \bigcup_{\omega\in\{1,\ldots,4\}^k} \psi_{\omega}(V))=0$. Therefore,
    \begin{align}
    \begin{aligned}\label{eq:inV}
        \vol_2\left({K}_{-\e}\cap V\right) 
        &= \sum_{k=0}^{\infty} \sum_{\omega\in\{1,\ldots,4\}^k} \vol_2\left({K}_{-\e}\cap \psi_{\omega}(W)\right)\\
        &= \sum_{k=0}^{\infty} \sum_{\omega\in\{1,\ldots,4\}^k} \vol_2\left((\psi_{\omega}{K})_{-\e}\cap \psi_{\omega}(W)\right)\\
        &= \sum_{k=0}^{\infty} \left(\frac49\right)^k \vol_2\left({K}_{-3^k\e}\cap W\right).
        \end{aligned}
    \end{align}
    Now, using the decomposition of $W$ into four congruent copies of $\Gamma$ and two congruent copies of $\Lambda$, see Fig.~\subref{fig:Gammatilde}, we obtain
    \begin{align}
    \begin{aligned}\label{eq:inW}
    &\vol_2\left({K}_{-3^k\e}\cap W\right)\\
    &\quad= 4\cdot \vol_2\left({K}_{-3^k\e}\cap \Gamma\right) +2\cdot \vol_2\left({K}_{-3^k\e}\cap \Lambda\right)\\
    &\quad= \begin{cases}
        4\cdot\frac{\sqrt{3}}{108}+2\cdot \frac{\sqrt{3}}{54} &\colon k\geq \lfloor-\frac{\log\e}{\log 3}\rfloor\\
        4\cdot\frac{\sqrt{3}}{108} + \frac{2\pi}3 3^{2k}\e^2 
        - 2\cdot 3^{2k}\e^2\arccos\left( \frac{1}{6\cdot 3^{k}\e}\right) 
        - \frac{\sqrt 3}{54}&\\
        \qquad
        +  \frac13\sqrt{3^{2k}\e^2-\frac{1}{36}}&\colon \lfloor-\frac12-\frac{\log\e}{\log 3}\rfloor\leq k< \lfloor-\frac{\log\e}{\log 3}\rfloor\\
        4\cdot \vol_2\left({K}_{-3^k\e}\cap \Gamma\right) +\frac{\pi}{3}3^{2k}\e^2 & \colon k< \lfloor-\frac12-\frac{\log\e}{\log 3}\rfloor.
    \end{cases}
    \end{aligned}
    \end{align}
    Combining \eqref{eq:inV} with \eqref{eq:inW} we can evaluate $\vol_2({K}_{-\e}\cap V)$, leading to \eqref{eq:innerkochvolume}. For this, note the following.
    \begin{enumerate}
        \item With $\alpha(\e):=\left\{ -\frac{\log{\e}}{\log{3}}\right\}$ as defined in the statement of this Lemma, it is convenient to write  $\lfloor-\frac{\log\e}{\log 3}\rfloor = -\frac{\log\e}{\log 3}-\alpha(\e)$.
        \item If $\alpha(\e)<\frac12$, then $\lfloor-\frac12-\frac{\log\e}{\log 3}\rfloor = \lfloor-\frac{\log\e}{\log 3}\rfloor-1$. If $\alpha(\e)\geq\frac12$, then $\lfloor-\frac12-\frac{\log\e}{\log 3}\rfloor = \lfloor-\frac{\log\e}{\log 3}\rfloor$.
        Thus, the middle case of \eqref{eq:inW} occurs if and only if $\alpha(\e)<\frac12$. 
        \item Due to the self-similarity of the Koch curve, $\vol_2 ({K}_{-\epsilon}\cap \Gamma) = 9\vol_2 ({K}_{-\epsilon/3}\cap \Gamma)$ whenever $\epsilon \leq \frac{1}{9}$.
    \end{enumerate}
\end{proof}
\begin{figure}
    \begin{subfigure}[t]{.48\textwidth}
        \centering
            \begin{tikzpicture}[scale=2.2,decoration=Koch snowflake]
    \draw[line width=0.01mm] decorate{decorate{decorate{ decorate{ decorate{ decorate{ (0,0) -- (3,0) }}}}}} -- ({1.5},{-sqrt(3)/2}) -- (0,0);
    \draw[gray,dashed] ({1/2},{-sqrt(3)/6}) -- (1,0);
    \draw[gray,dashed] ({5/2},{-sqrt(3)/6}) -- (2,0)  -- (2,{-sqrt(3)/3});
    \draw[gray,dashed] ({3/2},{-sqrt(3)/6}) -- (1,0) --  ({3/2},{sqrt(3)/6});
    \draw[gray,dashed] ({3/2},{-sqrt(3)/6}) -- (2,0) --  ({3/2},{sqrt(3)/6});
    \draw[gray,dashed] ({3/2},{-sqrt(3)/2}) -- ({3/2},{sqrt(3)/2});
    \draw (2,0) -- ({1.5},{-sqrt(3)/2}) node[midway,right]{$\frac13$};
    \draw (1,0) --  (1,{-sqrt(3)/3}) node[midway,right]{$\frac{\sqrt{3}}9$};
    \draw ({1/3},0) -- ({1.5/3},{-sqrt(3)/6}) node[midway,right]{$\frac19$};
    \end{tikzpicture}
    \caption{Visualisation of the lengths which lead to the different cases in the proof of Lem.~\ref{lem:innerKochflake}: 
    $\vol_2({K}_{-\e}\cap \Lambda)= \vol_2(\Lambda)$ if and only if $\e\geq\frac13$. 
    $\vol_2({K}_{-\e}\cap \Gamma)= \vol_2(\Gamma)$ if and only if $\e\geq\frac{\sqrt{3}}{9}$.
    $\vol_2({K}_{-\e}\cap \psi_1V)= \vol_2(\psi_1 V)$ if and only if $\e\geq\frac19$.}
    \label{fig:cases-proof}
    \end{subfigure}
    \hfill
\begin{subfigure}[t]{.48\textwidth}
    \centering
    \begin{tikzpicture}[scale=2.2,decoration=Koch snowflake]
        \draw[line width=0.01mm] decorate{decorate{decorate{ decorate{ decorate{ decorate{ (0,0) -- (3,0) }}}}}};
        \draw (0,0) -- ({1.5},{-sqrt(3)/2}) -- (3,0);
        \draw[gray,dashed] ({1/2},{-sqrt(3)/6}) -- (1,0) --  (1,{-sqrt(3)/3});
        \draw[gray,dashed] ({5/2},{-sqrt(3)/6}) -- (2,0)  -- (2,{-sqrt(3)/3});
        \draw[gray,dashed] ({3/2},{-sqrt(3)/6}) -- (1,0) --  ({3/2},{sqrt(3)/6});
        \draw[gray,dashed] ({3/2},{-sqrt(3)/6}) -- (2,0) --  ({3/2},{sqrt(3)/6});
        \draw[gray,dashed] ({3/2},{-sqrt(3)/2}) -- ({3/2},{sqrt(3)/2});
        \draw[fill=black] ([shift=(220:0.66cm)]2,0) arc (220:260:0.66cm) -- ({3/2},{-sqrt(3)/2}) -- ([shift=(280:0.66cm)]1,0) arc (280:320:0.66cm);
        \draw ([shift=(220:0.66cm)]2,0) -- (2,0) -- ([shift=(260:0.66cm)]2,0) ;
        \draw ([shift=(220:0.25cm)]2,0) arc (220:260:0.25cm);
        \node at (1.92,-0.15) {{\scriptsize{$\rho$}}};
        \node at (1.88,-0.4) {{\scriptsize{$\varepsilon$}}};
        \node at (1.2,-1/3) {$\Lambda$};
        \node at (0.8,-1/3) {$\Gamma$};
    \end{tikzpicture}
    \caption{Example of an inner $\e$ neighbourhood of ${K}$ for $\e\in(\frac{\sqrt{3}}{9}, \frac13]$ as in Case 2 in the proof of Lem.~\ref{lem:innerKochflake}. Here, $\psi_1V$ and $\Gamma$ are fully contained in ${K}_{-\e}$, whereas $\Lambda$ is not.}
    \label{fig:Case2}
    \end{subfigure}
\end{figure}
\begin{proof}[Proof of Thm.~\ref{thm:volume}.]
In this proof we fix $(k_1,k_2)\in\{0,\ldots,6\}^2$ and abbreviate the fractal spray $\Omega(k_1,k_2)$, the alphabet $\Sigma(k_1,k_2)$ and the IFS $\Phi(k_1,k_2)$ as introduced in Sec.~\ref{sec:kochspray} by $\Omega$, $\Sigma$ and $\Phi$ respectively. 
For $\beta\in [0,a)$, we define
\begin{equation*}
    N(a\ell+\beta):= \vol_2 \left( \Omega_{-\textup{e}^{-(a\ell+\beta)}}\right).
\end{equation*}
Recall from \eqref{eq:Omega} that $\Omega$ is a disjoint union of open sets
$\Omega = \bigcup_{k=0}^{\infty} \bigcup_{\omega\in\Sigma^k} \phi_{\omega}(G)$,
where the generator $G:=O\setminus\overline{\Phi O}$ can have several connected components, depending on $(k_1,k_2)$, and where $\Sigma^0:=\{\emptyset\}$ and $\phi_{\emptyset}$ is the identity on $\mathbb R^2$.
With the identity $\Omega_{-\e}:=\{ x\in\Omega\,:\,\dist(x,\partial\Omega)<\e\} = (\partial\Omega)_{\e}\cap \Omega$, where $F_{\e}:=\{ x\in \mathbb R^2\,:\, \dist(x,F)<\e \}$ denotes the \emph{$\e$-parallel set} of a set $F\subset\mathbb R^2$, we have that
\begin{align*}
    N(a\ell+\beta)
    &= \vol_2\Big{(} (\partial\Omega)_{\textup{e}^{-(a\ell+\beta)}}\cap \bigcup_{k=0}^{\infty}\bigcup_{\omega\in\Sigma^k} \phi_{\omega} G \Big{)}
    = \sum_{k=0}^{\infty} \sum_{\omega \in \Sigma^k} \vol_2 \left( (\partial\Omega)_{\textup{e}^{-(a\ell+\beta)}}\cap \phi_{\omega} G\right)\\
    &= \sum_{k=0}^{\infty} \sum_{\omega \in \Sigma^k} \vol_2 \left( (\phi_{\omega}(\partial\Omega))_{\textup{e}^{-(a\ell+\beta)}}\cap \phi_{\omega} G\right)
    =\sum_{k=0}^{\infty} \sum_{\omega \in \Sigma^k} \vol_2 \left( \phi_{\omega}\left((\partial\Omega)_{\textup{e}^{-(a\ell+\beta-a \nu_{\omega})}}\cap G\right)\right)\\
    &=\sum_{k=0}^{\infty} \sum_{\omega \in \Sigma^k} \textup{e}^{-2a\nu_{\omega}}\vol_2 \left( G_{-\textup{e}^{-(a(\ell-\nu_{\omega})+\beta)}}\right).
\end{align*}
In the last two equations \eqref{eq:nu} has been used.
Next, we consider the two-sided Fourier-Laplace transform $\widehat{N}_{\beta}$ dependent on $\beta\in[0,a)$, acting on $\mathbb C$ and given by
\begin{align*}
    \widehat{N}_{\beta}(z)
    &= \sum_{\ell=-\infty}^{\infty} \textup{e}^{a\ell z} N(a\ell +\beta).
\end{align*}
For $z\in\{z\in\mathbb C\mid 0<\Re(z)< 2-\dim_M(\partial\Omega)\}$, where $\dim_M(\partial\Omega)$ denotes the Minkowski dimension of $\partial\Omega$, the Fourier-Laplace transform $\widehat{N}_{\beta}(z)$ converges and the order of summation can be swapped, leading to the following conversion.
\begin{align*}
    \widehat{N}_{\beta}(z)
    &= \sum_{k=0}^{\infty} \sum_{\omega \in \Sigma^k} \textup{e}^{-2a\nu_{\omega}} \sum_{\ell=-\infty}^{\infty} \textup{e}^{a\ell z}\vol_2 \left( G_{-\textup{e}^{-(a(\ell-\nu_{\omega})+\beta)}}\right)\\
    &= \sum_{k=0}^{\infty} \sum_{\omega \in \Sigma^k} \textup{e}^{-a\nu_{\omega}(2-z)} \sum_{\ell=-\infty}^{\infty} \textup{e}^{a\ell z}\vol_2 \left( G_{-\textup{e}^{-(a\ell+\beta)}}\right)
\end{align*}
In the last equality we have used an index shift, as $\nu_{\omega}\in\mathbb Z$ by \eqref{eq:nu}.
Depending on $(k_1,k_2)$, the generator $\gen$ may have several connected components $K^{(0)},\ldots,K^{(k_1+k_2)}$, all of which are Koch snowflakes. Thus,
\begin{align*}
    \vol_2 \left( G_{-\textup{e}^{-(a\ell+\beta)}}\right)
    &= \sum_{j=0}^{k_1+k_2}\vol_2 \left( K^{(j)}_{-\textup{e}^{-(a\ell+\beta)}}\right)
    = \sum_{j=-k_1}^{k_2}b_j^2\vol_2 \left( K_{-\textup{e}^{-(a\ell+\beta+\log b_j)}}\right),
\end{align*}
with $b_{j}$ denoting the base length of the Koch snowflake $K^{(j)}$. In our setting, we have $b_0=1=\textup{e}^{-0\cdot a}$, $b_{j} = \sqrt{3}/3=\textup{e}^{-1\cdot a}$ for $j<0$ and $b_{j}=b_{k_1+k_2}=1/3=\textup{e}^{-2\cdot a}$ for $j>0$, implying
\begin{align*}
    \widehat{N}_{\beta}(z)
    &= (1+k_1\textup{e}^{a(z-2)}+k_2\textup{e}^{2a(z-2)})\sum_{k=0}^{\infty} \sum_{\omega \in \Sigma^k} \textup{e}^{-a\nu_{\omega}(2-z)} \sum_{\ell=-\infty}^{\infty} \underbrace{\textup{e}^{a\ell z}\vol_2 \left( K_{-\textup{e}^{-(a\ell+\beta)}}\right)}_{=:h(\ell)}.
\end{align*}
We can evaluate the series with indices $k$ and $\ell$ independently. For the series with index $k$ we use that $\textup{e}^{-a\nu_{\omega}}=r_{\omega}$ and that $\sum_{\omega\in\Sigma^k}r_{\omega}^{2-z}= \big{(} \sum_{i\in\Sigma} r_i^{2-z} \big{)}^k$. 
For the series with index $\ell$, we use \eqref{eq:innerkochvolume}, and split the series in the following way. 
$\sum_{\ell=-\infty}^{\infty}h(\ell)
= \sum_{\ell=-\infty}^{1}h(\ell)
+ h(2) + h(3)
+\sum_{\ell=4}^{\infty}h(\ell)\mathds{1}_{2\mathbb Z}(\ell)
+\sum_{\ell=4}^{\infty}h(\ell)\mathds{1}_{2\mathbb Z+1}(\ell)$.
\begin{align*}
    &\widehat{N}_{\beta}(z) \cdot (1+k_1\textup{e}^{a(z-2)}+k_2\textup{e}^{2a(z-2)})^{-1}\\
    &\quad= \frac{1}{1-\sum_{i \in \Sigma} r_i^{2-z}} 
    \bigg[
    \frac{\textup{e}^{a z}}{1-\textup{e}^{-az}}\cdot \frac{2\sqrt{3}}{5} 
    + \frac{\textup{e}^{2az}}{3}\cdot\left( \frac{7\sqrt{3}}{10} + \sqrt{\textup{e}^{-2\beta}-\frac{1}{4}} + 2\textup{e}^{-2\beta}\arcsin{\left(\frac{\textup{e}^{\beta}}{2}\right)} - \frac{\pi \textup{e}^{-2\beta}}{3} \right)\\
    &\quad\qquad\qquad\qquad \qquad
    + \textup{e}^{3az}\cdot\left( \frac{8\sqrt{3}}{45} + \frac{\pi \textup{e}^{-2\beta}}{27} + 12\vol_2(K_{-\textup{e}^{-\beta}\sqrt{3}^{-3}}\cap \Gamma)\right)\\
    &\quad\qquad\qquad\qquad \qquad
    + u\left(\frac{\beta}{2a}\right)\textup{e}^{-\beta(2-\frac{\log 4}{\log 3})}\cdot \frac{\textup{e}^{4a(z-2+\frac{\log 4}{\log 3})}}{1 - \textup{e}^{2a(z-2+\frac{\log 4}{\log 3})}} + v\left(\frac{\beta}{2a}\right) \textup{e}^{-2\beta}\frac{\textup{e}^{4a (z-2)}}{1-\textup{e}^{2a(z-2)}} \\
    &\quad\qquad\qquad\qquad \qquad
    + \widetilde{u}\left(\frac{a+\beta}{2a}\right)\textup{e}^{-\beta(2-\frac{\log 4}{\log 3})}\cdot \frac{\textup{e}^{5a(z-2+\frac{\log 4}{\log 3})}}{1 - \textup{e}^{2a(z-2+\frac{\log 4}{\log 3})}} + v\left(\frac{a+\beta}{2a}\right) \textup{e}^{-2\beta}\frac{\textup{e}^{5a (z-2)}}{1-\textup{e}^{2a(z-2)}} 
    \bigg]\\
    &\quad=: \frac{1}{1-\sum_{i \in \Sigma} r_i^{2-z}} \cdot L(z)
\end{align*}
The right-hand side has a meromorphic extension to $\mathbb C$ with simple poles at $z$ in
\begin{equation*}
\mathcal Z :=\bigg{\{} z\in\mathbb C \mid \sum_{i \in \Sigma} r_i^{2-z} =1 \bigg{\}} \quad \cup \quad \mathcal S:=\bigg{\{} 0, 2 -\frac{\log 4}{\log 3}, 2\bigg{\}}.
\end{equation*}
Define
\begin{align}
\begin{aligned}\label{eq:Minkcoeff}
    {\coeffmink}_{\beta}(2)&:= \frac{\textup{e}^{-2\beta}}{2(1-\#\Sigma)}\cdot \left[ v\left(\frac{\beta}{2a}\right) + v\left(\frac{a+\beta}{2a}\right)\right],\\
    {\coeffmink}_{\beta}\left(2-\frac{\log 4}{\log 3}\right)
    &:= \frac{(1+\frac{k_1}{2}+\frac{k_2}{4})\textup{e}^{-\beta(2-\frac{\log 4}{\log 3})}}{2\big{(}1-\sum_{i\in\Sigma}r_i^{{\log 4}/{\log 3}}\big{)}}\cdot \left[ u\left(\frac{\beta}{2a}\right) + \widetilde{u}\left(\frac{a+\beta}{2a}\right) \right]\quad\text{and}\\[1.2em]
    {\coeffmink}_{\beta}(z)
    &:= -\frac{a(1+k_1\textup{e}^{a(z-2)}+k_2\textup{e}^{2a(z-2)})}{\sum_{i\in \Sigma}\log r_i \cdot r_i^{2-z}}\cdot L(z)
    \end{aligned}
\end{align}
for $z\in\mathcal Z$.
Then 
\begin{equation*}
    H_{\beta} (s) 
    := \widehat{N}_{\beta}(s) - \sum_{\ell =-\infty}^{-1} \textup{e}^{a\ell s} N(a\ell +\beta) - \frac{{\coeffmink}_{\beta}(2)}{1-\textup{e}^{a(s-2)}} - \frac{{\coeffmink}_{\beta}(2- \frac{\log 4}{\log 3})}{1-\textup{e}^{a(s-2+\frac{\log 4}{\log 3})}} - \sum_{z\in \mathcal{Z}_{\textup P}}\frac{{\coeffmink}_{\beta}(z)}{1-\textup{e}^{a(s-z)}}
\end{equation*}
extends to a holomorphic function on $\mathbb C$, where $\mathcal{Z}_{\textup P} := \{ z \in \mathcal{Z}: \Im(z) \in [0,\frac{2\pi}{a}) \}$. On $\{z\in\mathbb C\mid  0<\Re(z)< \min_{s\in\mathcal Z} \Re(s) \}$ each summand of $H_{\beta}$ can be developed into a power series:
\begin{equation*}
    H_{\beta}(s)
    = \sum_{\ell=0}^{\infty}\textup{e}^{a\ell s} \left[ N(a\ell+\beta) - {\coeffmink}_{\beta}(2)\textup{e}^{-2a\ell} - {\coeffmink}_{\beta}\left(2- \frac{\log 4}{\log 3}\right)\textup{e}^{-a\ell(2-\frac{\log 4}{\log 3})} - \sum_{z\in \mathcal{Z}_{\textup P}} \textup{e}^{-a\ell z}{\coeffmink}_{\beta}(z) \right]
\end{equation*}
As a power series has a singularity on its radius of convergence, $H_{\beta}$ being holomorphic on $\mathbb C$ implies that as $\ell\to\infty$
\begin{equation*}
     N(a\ell+\beta) 
     = {\coeffmink}_{\beta}(2)\textup{e}^{-2a\ell} 
     + {\coeffmink}_{\beta}\left(2- \frac{\log 4}{\log 3}\right)\textup{e}^{-a\ell(2-\frac{\log 4}{\log 3})} 
     +\sum_{z\in \mathcal{Z}_{\textup P}} \textup{e}^{-a\ell z}{\coeffmink}_{\beta}(z)
     +o(\textup{e}^{-a\ell\gamma})
\end{equation*}
for any $\gamma>0$.

For a selection of choices of $(k_1,k_2)$, the coefficients ${\coeffmink}_\beta (z)$ at $z \in \mathcal{S}$ take the values as shown in Tab.~\ref{tab:valuecontentG}.
\begin{table}[h]
 \centering
 \begin{tabular}{c|c|c}
 $(k_1,k_2)$ & Values of ${\coeffmink}_\beta(2)$ & Values of ${\coeffmink}_\beta(2-\frac{\log4}{\log3})$\\
 \hline
 $(0,0)$ & $-\frac{\textup{e}^{-2\beta}}{22}\cdot\left[ v\left(\frac{\beta}{2a}\right) + v\left(\frac{a+\beta}{2a}\right)\right]$ & $-\frac{2\textup{e}^{-\beta(2-\frac{\log 4}{\log 3})}}{5}\cdot\left[ u\left(\frac{\beta}{2a}\right) + \widetilde{u}\left(\frac{a+\beta}{2a}\right) \right]$ \\
 $(0,6)$ & $-\frac{\textup{e}^{-2\beta}}{82}\cdot\left[ v\left(\frac{\beta}{2a}\right) + v\left(\frac{a+\beta}{2a}\right)\right]$ & $-\frac{10\textup{e}^{-\beta(2-\frac{\log 4}{\log 3})}}{13}\cdot\left[ u\left(\frac{\beta}{2a}\right) + \widetilde{u}\left(\frac{a+\beta}{2a}\right) \right]$\\
 $(6,6)$ & $-\frac{\textup{e}^{-2\beta}}{142}\cdot\left[ v\left(\frac{\beta}{2a}\right) + v\left(\frac{a+\beta}{2a}\right)\right]$ & $-\frac{22\textup{e}^{-\beta(2-\frac{\log 4}{\log 3})}}{19}\cdot\left[ u\left(\frac{\beta}{2a}\right) + \widetilde{u}\left(\frac{a+\beta}{2a}\right)\right]$
\end{tabular}
\caption{Exact values of the coefficients of the asymptotic expansion of the Lebesgue measure of an inner tubular neighbourhood of $\Omega(k_1,k_2)$.}\label{tab:valuecontentG}
\end{table}
\end{proof}

 \bibliographystyle{amsplain}
 \providecommand{\bysame}{\leavevmode\hbox to3em{\hrulefill}\thinspace}
\providecommand{\MR}{\relax\ifhmode\unskip\space\fi MR }

\end{document}